\documentclass{article}

\usepackage{amsmath,amsthm,amssymb} 
\usepackage{enumerate,epic,eepic}

\newtheorem{thm}{Theorem}
\newtheorem{prop}[thm]{Proposition}
\newtheorem{lemma}[thm]{Lemma}

\allowdisplaybreaks[2]

\DeclareFontEncoding{OT2}{}{}
\DeclareFontSubstitution{OT2}{wncyss}{m}{n}
\DeclareSymbolFont{cyss}{OT2}{wncyss}{m}{n}
\DeclareMathSymbol{\sh}{\mathbin}{cyss}{`x}

\title{The Sum Formula of Multiple Zeta Values and Connection Problem of the 
Formal Knizhnik-Zamolodchikov Equation}
\author{Jun-ichi OKUDA and Kimio UENO}
\date{February 17, 2004}

\begin{document}
\maketitle
\begin{abstract}
The sum formula for multiple zeta values are derived, via the Mellin transform, 
from the Euler connection formula and the Landen connection formula for 
polylogarithms. These connection formulas for multiple polylogarithms
will be considered in the framework of 
the theory of the formal Knizhnik-Zamolodchikov equation.
\end{abstract}

\section{Introduction}

In this paper, we will derive the sum formula \cite{G} for multiple zeta values (MZVs, for
short) in two ways: We will show that the Mellin transform of the Euler connection
formula and the Landen connection formula for polylogarithms 
give the sum formula, respectively. Furthermore we will clarify the meaning of 
these connection formulas for multiple polylogarithms (MPLs, briefly) in the framework
of the connection problem of the formal Knizhnik-Zamolodchikov equation.

In \cite{OU}, we considered the Ohno relation for MZVs \cite{O} by means of generating functional
method: That is, we introduced the two generating functions $f((a_i,b_i)_{i=1}^s;\lambda)$,
$g((a_i,b_i)_{i=1}^s;\lambda)$ by which the Ohno relation is represented as $f=g$, and
found the subfamily $F_{\mathbf{k}}(\lambda)=G_{\mathbf{k}}(\lambda)$ of the Ohno relation
which is converted to the Landen connection formula for MPLs by the inverse 
Mellin transform. This subfamily was called the reduced Ohno relation.

The multiple polylogarithm is, by definition,
\begin{align}\label{MPL}
  \mathrm{Li}_{k_1,k_2,\dotsc,k_m}(z) = \sum_{n_1>n_2>\dotsb>n_m>0}
                                \frac{z^{n_1}}{n_1^{k_1}n_2^{k_2}\dotsm n_m^{k_m}},
\end{align}
and the Landen connection formula \cite{OU}
\begin{align}\label{LANDENMPL}
  \mathrm{Li}_{k_1,\ldots,k_m}(z) = (-1)^n \sum_{\substack{c_1,\dotsc,c_m\\ \text{weight of} \, c_j=k_j}}
                                 \mathrm{Li}_{c_1,\dotsc,c_m} \left(\frac{z}{z-1}\right)
\end{align}
can be thought of the connection formula between $z=1$ and $z=\infty$ for MPLs.
On the other hand, the Euler inversion formula \cite{L} for the dilogarithm
\begin{align}\label{EULERDI}
     \mathrm{Li}_2(z)+\mathrm{Li}_1(1-z)\mathrm{Li}_1(z)+\mathrm{Li}_2(1-z)=\zeta(2)
\end{align}
and its generalization (the Euler connection formula) give the connection 
between $z=0$ and $z=1$ for MPLs. Only the polylogarithm case will be treated
in this paper. It seems interesting to calculate relations for general MPLs. 
 
 These connection formulas can be easily verified by 
elementary methods, but can be understood in a unified way by considering the connection problem
for the formal Knizhnik-Zamolodchikov equation 
\begin{align}\label{FKZ}
   \frac{dG}{dz} = \left( \frac{X}{z}+\frac{Y}{1-z} \right) G. \tag{KZ}
\end{align}

This paper is organized as follows: In Section 2, we will briefly review the theory of
the shuffle algebra $\mathfrak{h}=(\mathbb{C}\langle x,y \rangle, \sh )$
and the regularization map defined on it, according to the arguments in \cite{IK, mph} 
and \cite{okuda}. As an element in
$\mathfrak{H}=\mathbb{C}\langle\langle X,Y \rangle\rangle$, which is 
thought of to be the topological dual of $\mathfrak{h}$,
the Drinfel$'$d associator $\varphi_{\mathrm{KZ}}$ constructed in \cite{D}
is canonically introduced.
In Section 3, we will consider the link between the connection formulas for MPLs and
the connection problem for \eqref{FKZ}. 
The solutions $G_0$, $G_1$, $G_{\infty}$ to \eqref{FKZ} are defined to be 
unique ones satisfying the prescribed asymptotic properties around each singular point.
In particular, $G_0$ plays a role of the generating function for MPLs, and 
$G_1$ and $G_{\infty}$ are expressed in terms of $G_0$ and homomorphisms induced from 
linear fractional transformations preserving the singular points. Thus the connection problem
for \eqref{FKZ} leads to the connection formulas for MPLs. Furthermore, we can verify 
the so called hexagon relation \cite{D} satisfied by the Drinfel$'$d associator
by using a sort of the braid relation for the induced homomorphisms.
These results are parallel to those in \cite{mph}, however they did not use the 
regularization map defined on $\mathfrak{h}$ which was firstly introduced by \cite{IK} and
\cite{Ka}.
In Section 4, we will see that the sum formula for MZVs follows from the Mellin 
transform of the Euler connection formula and the Landen connection formula for polylogarithms.

There are two appendices: In Appendix 1, we will consider the connection problem 
for the system ($\mathrm{KZ}_3$) \cite{D} over the configuration space 
of ordered three points in the complex line \ 
$\mathbb{F}_3(\mathbb{C})=\{(z_1,z_2,z_3) \in \mathbb{C}^3 \,|\, z_i \neq z_j \ (i \neq j) \, \}$
\begin{align}\label{KZ3}
   dW = \left( \sum_{1\leq i<j \leq 3}X_{ij}d\log(z_i-z_j) \right)W,  \tag{KZ$_3$}
\end{align}
where the coefficients satisfy $X_{ij}=X_{ji}, \quad [X_{ij}, X_{ik}+X_{kj}]=0$,
with the aid of the connection problem for \eqref{FKZ}. Through these investigations, we see that
the relations for MPLs which come from the connection problem of \eqref{KZ3} are only
the Euler connection formula and the Landen connection formula in a essential sense.
In Appendix 2, written by the second author and Michitomo Nishizawa,
we will give an alternative proof for the functional relation of
the Hurwitz zeta function by using the connection problem for
the confluent hypergeometric equation. 
This proof appeared first in \cite{UN1} in an implicit form, and second in \cite{UN2} 
which was written unfortunately in Japanese. In these articles we considered mainly
the functional relation for the q-analogue of the Hurwitz zeta function.
So we decided to publish again the proof in a complete form. \\

\noindent
\textbf{Acknowledgement} \quad The authors would like to express their deep gratitude
to the organizing committee of the conference ZTQ, especially to Professor Yasuo Ohno,
for giving the authors the opportunity presenting their research.  

The second author is
partially supported by JPSP Grant-in-Aid No. 15540050 and
Waseda University Grant for Special Research Project (2002A-067,
2003A-069).

\section{Shuffle Algebra}

We define the shuffle product $\sh$
on the non-commutative  polynomial algebra $\mathfrak{h}=\mathbb{C} \langle x, y \rangle$
of letters $x$ and $y$, 
recursively by
\begin{align*}
 1 \sh w & = w \sh 1 = w,\\
 l_1 w_1 \sh l_2 w_2
 & = l_1 (w_1 \sh l_2 w_2) + l_2 (l_1 w_1 \sh w_2),
\end{align*}
where $w, w_1$ and $w_2$ are monomials
and $l_1$ and $l_2$ are the letters.
Now let us introduce the subalgebras $\mathfrak{h} \supset\mathfrak{h}^1\supset\mathfrak{h}^0$
with respect to \ $\sh$ and concatenation product
defined by
\begin{align}
 \mathfrak{h}^1 = \mathbb{C}1 \oplus \mathfrak{h} y, \qquad  
 \mathfrak{h}^0 = \mathbb{C}1 \oplus x \mathfrak{h} y.
\end{align}
For $w=x^{k_1-1}y\dotsm x^{k_m-1}y\in\mathfrak{h}^0$,
we set a linear map $\widehat{\zeta}:\mathfrak{h}^0\to\mathbb{C}$ by
\begin{align}
\widehat{\zeta}(w)=\zeta(k_1,\dotsc,k_m)= \sum_{n_1>\dotsb>n_m>0}
                  \frac{1}{n_1^{k_1}\dotsm n_m^{k_m}},
\end{align}
which is a multiple zeta value.
Then $\widehat{\zeta}$ is a homomorphism with respect to $\sh$,
that is,
\begin{align*}
 \widehat{\zeta}(w_1 \sh w_2)
 = \widehat{\zeta}(w_1) \widehat{\zeta}(w_2).
 \label{eq:sh-hom}
\end{align*}

It is well known that $\mathfrak{h}$ is a polynomial algebra
with respect to the products $\sh$ and $\star$
generated by the Lyndon words (cf.\ \cite{reut}).
In particular the Lyndon words contain the letters $x$ and $y$, and
$\mathfrak{h}$ can be decomposed as follows:
\begin{align}
 \mathfrak{h} &= \bigoplus_{n=0}^\infty \mathfrak{h}^1 \sh x^{\sh n}
    = \mathfrak{h}^1[x]\\
    &= \bigoplus_{m,n=0}^\infty y^{\sh m} \sh \mathfrak{h}^0 \sh x^{\sh n}
    = \mathfrak{h}^0[x,y].
    \label{eq:h decomp}
\end{align}
Let $\mathrm{reg} : \mathfrak{h}=\mathfrak{h}^0[x,y]\to\mathfrak{h}^0$ be 
a map to associate the constant term in the decomposition $\mathfrak{h}=\mathfrak{h}^0[x,y]$
to each element in $\mathfrak{h}$. This is called the reguralization map,
and is a $\sh$-homomorphism.
\begin{prop}[\cite{IK}]
For $w \in \mathfrak{h}^0$, 
\begin{align}
 \mathrm{reg}(y^m w x^n)
  &= \sum_{i=0}^m \sum_{j=0}^n (-1)^{i+j}
   y^{i} \sh y^{m-i} w x^{n-j} \sh x^j,\\
 \intertext{or equivalently}
 y^m w x^n
  &= \sum_{i=0}^m \sum_{j=0}^n
   y^{i} \sh \mathrm{reg} (y^{m-i} w x^{n-j}) \sh x^j. \label{CHEN}
 \end{align}
\end{prop}

Let us introduce the algebra $\mathfrak{H} = \mathbb{C}\langle\langle X, Y\rangle\rangle$
of non-commutative formal power series over $\mathbb{C}$ in the letters $X$ and $Y$.
This is thought of to be the topological dual of $\mathfrak{h}$.
Then the latter equation (\ref{CHEN}) is expressed by the ``formal'' Chen series 
in $\mathfrak{h} \widehat{\otimes} \mathfrak{H}$ \cite{IK}
\begin{align}\label{ELEMENT}
 \sum_W wW &= 1 + xX + yY + xx XX + xy XY + yx YX + yy YY + \dotsb \nonumber \\
 &=
 \exp_\sh (x X) \left(\sum_W \mathrm{reg}(w)W\right)  \exp_\sh (yY),
\end{align}
where $W$ runs over all the monomials in $\mathfrak{H}$,
and $\exp_\sh$ is a series defined by
\begin{align*}
 \exp_\sh(xX)
  &= 1 + xX + \frac{x^{\sh 2} X^2}{2!} + \frac{x^{\sh 3} X^3}{3!} + \dotsb\\
  &= 1 + xX + x^2 X^2 + x^3X^3 + \dotsb
\end{align*}
as well as $\exp_\sh(yY)$.

Let $\varphi_{\mathrm{KZ}}(X,Y)\in\mathfrak{H}$ be the series \cite{IK}
\begin{align}\label{KZASS}
 \varphi_{\mathrm{KZ}}(X,Y) := \sum_W \widehat{\zeta}(\mathrm{reg}(w)) W.
 \end{align}
According to \cite{Kas}, we refer to $\varphi_{\mathrm{KZ}}$ as the Drinfel$'$d associator.

We define the antipode $S$ of $\mathfrak{h}$ and $\mathfrak{H}$ to be an 
anti-involution such that
\begin{align*}
 S: x&\mapsto -x, \qquad y \mapsto -y,\\
    X&\mapsto -X, \qquad Y \mapsto -Y.
\end{align*}
We note that $S$ is a $\sh$-homomorphism on $\mathfrak{h}$ and yields
the inverse of the generating series:
\begin{align}
 \left(\sum_W wW\right)^{-1} = \sum_W S(w)W = \sum_W wS(W).
\end{align}

\section{Multiple Polylogarithms and the formal KZ equation}

For $w=x^{k_1-1}yx^{k_2-1}y \dotsm x^{k_m-1}y \in \mathfrak{h}^1$ and $x$ 
we define the multiple polylogarithms
on $\mathbb{C} \backslash (-\infty,0]\cup[1,\infty)$ by
\begin{align}
 \mathrm{Li}(x;z) &= \log z,\\
 \mathrm{Li}(w;z)
  & = \mathrm{Li}_{k_1,k_2,\dotsc,k_m}(z) \nonumber \\
  & = \sum_{n_1>n_2>\cdots>n_r>0}
    \frac{z^{n_1}}{n_1^{k_1}n_2^{k_2}\dotsm n_m^{k_m}}
     \qquad (|z|<1)   \nonumber \\
  &=
   \underbrace{\int_0^z\frac{dz}{z} \dotsb \int_0^z\frac{dz}{z}}_{k_1-1}
    \int_0^{z}\frac{dz}{1-z}
   \dotsb
   \underbrace{\int_0^z\frac{dz}{z} \dotsb \int_0^z\frac{dz}{z}}_{k_r-1}
    \left( -\log(1-z) \right),
\end{align}
where the branch of $\log z$ is chosen as the principal value.
Because
$\mathrm{Li}(\bullet;z)$ on $\mathfrak{h}^1$ is a $\sh$-homomorphism \cite{ree}
and $\mathfrak{h}$ can be decomposed as \eqref{eq:h decomp},
$\mathrm{Li}(\bullet;z)$ extends from $\mathfrak{h}^1$ to the whole $\mathfrak{h}$
as a $\sh$-homomorphism.
For word $w\in\mathfrak{h}^0$
the evaluation of $\mathrm{Li}(w;z)$ tending $z\to 1-0$ is
\begin{align*}
 \lim_{z\to 1-0} \mathrm{Li}(w;z) = \widehat{\zeta}(w).
\end{align*}
From this evaluation, \eqref{eq:sh-hom} obviously follows.

Applying $\mathrm{Li}(\bullet;z)$ to the Chen series $\sum_W wW$,
we have a $\mathfrak{H}$-valued function
\begin{align}\label{GO}
 G_0(X,Y;z) = G_0(z) = \mathrm{Li}(\bullet;z)\left( \sum_W wW \right)
  = \sum_W \mathrm{Li}(w;z) W
 \end{align}
 holomorphic on $\mathbb{C} \backslash (-\infty,0]\cup[1,\infty)$. 
 By virtue of \eqref{ELEMENT}, we have
 \begin{align}\label{GOREG}
 G_0(X,Y;z) &=
  \exp_{\sh}\left( \mathrm{Li}(y;z)Y \right)
   \left( \sum_W \mathrm{Li}(\mathrm{reg}(w);z)W \right)
    \exp_{\sh}\left( \mathrm{Li}(x;z)X \right) \nonumber \\
 &=
  \exp\left( -Y\log(1-z) \right)
   \left( \sum_W \mathrm{Li}(\mathrm{reg}(w);z)W \right)
    \exp\left( X\log z \right) \nonumber \\
 &=
  (1-z)^{-Y} \overline{G}_0(X,Y;z) z^X 
\end{align}
where
\begin{align*}
 \overline{G}_0(X,Y;z) = \sum_W \mathrm{Li}(\mathrm{reg}(w);z)W.
\end{align*}
$G_0(X,Y;z)$ is a unique solution to the system \eqref{FKZ}
\begin{align*}
   dG = \Omega_{\mathrm{KZ}} G,
 \tag{KZ}
\end{align*}
where
\begin{align*}
 \Omega_{\mathrm{KZ}} = X d\log z - Y d\log(1-z),
\end{align*}
with the asymptotic behavior
\begin{align}
 G_0(X,Y;z) \times z^{-X}  \longrightarrow 1 \qquad (z \to 0)
\end{align}
(cf.\ \cite{mph, okuda}).
There also exist unique solutions $G_1$ and $G_{\infty}$ of \eqref{FKZ}
with the following asymptotic behavior
\begin{alignat}{4}
 G_1(X,Y;z) &\times (1-z)^{Y}      &\longrightarrow  1& \qquad&(z \to 1),\\
 G_{\infty}(X,Y;z)& \times (1/z)^{-X+Y}  &\longrightarrow  1& \qquad&(z \to \infty).
\end{alignat}
Let us introduce the connection matrices $C^{(jk)}$ by
\begin{align}
 G_k = G_j \times C^{(jk)}, \qquad \text{($j,k=0,1$ or $\infty$)}.
\end{align}
We investigate the explicit form of the solutions $G_1$, $G_{\infty}$, 
and the connection matrices (the connection problem for \eqref{FKZ}!).
To do this, we use the linear fractional transformations
following \cite{mph, okuda}.
There are six functions preserving the singular points of \eqref{FKZ}:
\begin{align}
 &z,& &1-z,& &\frac{1}{z},&
 &\frac{z}{z-1},& &\frac{1}{1-z},& &\frac{z-1}{z}.\\
\intertext{Each function corresponds to the permutation of singular points:}
 &e,& &(01),& &(0\infty),&
 &(1\infty),& &(01\infty),& &(10\infty).
\end{align}
The linear fractional transformations $f$ acts on functions by the pull-back
\begin{align}
  (f^*G)(z) = G(f(z)).
\end{align}
We also define the actions on $\mathfrak{h}$,
$\mathfrak{H}$, which is denoted by $f^*$, $f_*$ respectively, by the following rule: 
Let us identify the letter $x$ as $d\log z$, and the letter $y$ as $-d\log(1-z)$. 
Then $\Omega_{\mathrm{KZ}}= xX + yY$. 
We define $f^*(x)$ and $f^*(y)$ to be the pull-back $f^*d\log z = d\log f(z)$, 
$-f^*d\log(1-z)=-d\log(1-f(z))$, respectively. Note that these one forms are again
expressed as linear combinations of $x$ and $y$. We write it as
\begin{align*}
   (f^*(x),f^*(y)) = (x,y)A(f),
\end{align*}
where $A(f)$ is a numerical $2 \times 2$ matrix. Then we define $f_*(X)$ and $f_*(Y)$ as
\begin{align*}
      \left(
        \begin{array}{c}
            f_*(X) \\
            f_*(Y)
      \end{array}
           \right)
           = A(f)
           \left(
        \begin{array}{c}
            X \\
            Y
      \end{array}
           \right).
\end{align*}
 So we have
\begin{align}
      f^*\Omega_{\mathrm{KZ}} = f^*(x) X + f^*(y) Y = x f_*(X) + y f_*(Y).
\end{align}
Then we have the following table:
\begin{align}
 \begin{array}[tb]{c||c|c|c|c}
               & x    & y    & X    & Y \\
  \hline
  \hline
 z             & x    & y    & X    & Y \\
  \hline
 1-z           & -y   & -x   & -Y   & -X \\
  \hline
 \frac{1}{z}   & -x   & x+y  & -X+Y & Y \\
  \hline
 \frac{z}{z-1} & x+y  & -y   & X    & X-Y \\
  \hline
 \frac{1}{1-z} & y    & -x-y & -Y & X-Y\\
 \hline
 \frac{z-1}{z} & -x-y & x    & -X+Y  & -X,
  \end{array}
\end{align}
e.g.\ 
\begin{align*}
  1-z:
  x&\mapsto -y,& y&\mapsto -x,&
  X&\mapsto -Y,& Y&\mapsto -X.&
\end{align*}
These extends homomorphisms on $\mathfrak{h}$, and $\mathfrak{H}$ respectively, and satisfy
\begin{align}
 & (f \circ g)^* = g^* \circ f^*, \qquad (f \circ g)_* = f_* \circ g_* \\
 & \sum_W f^*(w) \ W = \sum_W w \ f_*(W)
\end{align}
and $f^*$ is a $\sh$-homomorphism of $\mathfrak{h}$. 
Note that $f_* \circ g^* = g^* \circ f_*$ for any $f, g$.

For any $\mathfrak{H}$-valued function $F(X,Y,;z)=\sum_{W} \ell (w;z) W$, we set
\begin{align*}
 \pi(f)F(X,Y,;z) &= f_* (f^*)^{-1}F(X,Y;z)\\
              &= F(f_*(X),f_*(Y);f^{-1}(z))\\
              &= \sum_{W} \ell (w;f^{-1}(z))f_*( W)\\
              &= \sum_{W} \ell (f^*(w);f^{-1}(z))W.
\end{align*}
Since $H = (f^*)^{-1} G_0$ is a solution of the equation $dH = (f^*)^{-1}\Omega_{\mathrm{KZ}} \ H$,
applying $f_*$ to $(f^*)^{-1} G_0$,
\begin{align}
\pi(f) G_0(X,Y;z) = \sum_W \mathrm{Li}(f^*(w);f^{-1}(z))W
\end{align}
becomes a solution of \eqref{FKZ} again. Note that this transformation is $\mathbb{C}$-linear,
but not $\mathfrak{H}$-linear, and satisfy
\begin{align}
      \pi(f \circ g) = \pi(f) \circ \pi(g)
\end{align}
for any $f$ and $g$.
The asymptotic behavior of this solution around $f(0)$ can be obtained by 
\begin{align*}
& \pi(f) G_0(X,Y;z)\\
  & =
 \left(1-f^{-1}(z)\right)^{-f_*(Y)}
  \   \overline{G}_0(f_*(X), f_*(Y); f^{-1}(z))
   \left(f^{-1}(z)\right)^{f_*(X)}.
\end{align*}
Thus, from the uniqueness of the solutions satisfying the prescribed asymptotic behavior,
\begin{prop} We have 
\begin{align}
 \pi(1-z)(G_0) &= G_1,  \\
 \pi\left(\frac{1}{z}\right)(G_0) &= G_\infty,\\
 \pi\left(\frac{z}{z-1}\right)(G_0) &= G_0 \times \exp(\mp X\pi i),
  \label{eq:full landen}\\
 \pi\left(\frac{1}{1-z}\right)(G_0) &= G_1 \times \exp(\mp Y\pi i),\\
 \pi\left(\frac{z-1}{z}\right)(G_0)
   &= G_\infty \times \exp(\pm (-X+Y)\pi i),
\end{align}
where $z\in\mathbf{H}_{\pm}$.
\end{prop}
We should observe that in the latter three cases
all the cut of the logarithms which appear in the singular parts of the solutions
separate $\mathbb{C}$ to $\mathbf{H}_{+}\cup\mathbf{H}_{-}$, for example, 
\begin{align*}
 \log\left(\frac{z}{z-1}\right)+\log(1-z)-\log z = \mp \pi i  \quad (z \in \mathbf{H}_{\pm}).
\end{align*}
In particular \eqref{eq:full landen} is explicitly written as
\begin{align}
 \sum_W \mathrm{Li}\left(\left(\frac{z}{z-1}\right)^*(w);\frac{z}{z-1} \right) W
 =
 \left(\sum_W \mathrm{Li}(w;z) W \right) \times \exp(\pm X\pi i),
\end{align}
where
\begin{align*}
 \left(\frac{z}{z-1}\right)^*&:& x &\mapsto x+y,& y &\mapsto -y,
\end{align*}
so that the equality of the coefficient of $W \ (\text{a word} \in \mathfrak{H} Y)$ 
gives rise to the Landen connection formula \eqref{LANDENMPL}:
\begin{prop} We have
\begin{align}
 \mathrm{Li}\left(\left(\frac{z}{z-1}\right)^*(w);\frac{z}{z-1}\right)= \mathrm{Li}(w;z),
\end{align}
where $w$ is a word $\in \mathfrak{h}y$. 
\end{prop}
Using the antipode $S$, we can describe the inverse of $G_1$ as follows:
\begin{align*}
 & G_1(X,Y;z)^{-1}
 =
 \sum_W \mathrm{Li}\left((1-z)^*(w);1-z \right) S(W)\\
 &=
 \sum_W \mathrm{Li}\left( (1-z)^*\circ S(w);1-z \right) W\\
 &=
 \exp\left(-X \log z \right)
  \left(
   \sum_W \mathrm{Li}\left( \mathrm{reg}\circ(1-z)^*\circ S(w);1-z \right) W
  \right)
    \exp\left(Y \log (1-z) \right)\\
 &=
 \exp\left(-X \log z \right)
  \overline{G}_0(-Y,-X;1-z)^{-1}
    \exp\left(Y \log (1-z) \right).
\end{align*}
so that the ratio of these two solutions is computed, by letting $z \to 0$ and $z \to 1$, 
as follows:
\begin{align}
 C^{(10)} &= G_1(X,Y;z)^{-1} G_0(X,Y;z)
  \label{eq:full euler}\\
 &=
  \exp\left(-X \log z \right)
   \overline{G}_0(-Y,-X;1-z)^{-1}
    \overline{G}_0(X,Y;z)
     \exp\left(X \log z \right) \nonumber \\
 &=
  \overline{G}_0(X,Y;1)=\overline{G}_0(-Y,-X;1)^{-1}. \nonumber
\end{align}
Hence we obtain 
\begin{align*}
 C^{(10)} = \overline{G}_0(X,Y;1) = \sum_W\mathrm{Li}(\mathrm{reg}(w);1)W
                                  = \sum_W\widehat{\zeta}(\mathrm{reg}(w))W
                                  = \varphi_{\mathrm{KZ}}(X,Y)
\end{align*}
and
\begin{align*}
  C^{(10)} &= \overline{G}_0(-Y,-X;1)^{-1} =  \varphi_{\mathrm{KZ}}(-Y,-X)^{-1} \\
           &= \sum_W\mathrm{Li}(\mathrm{reg}\circ\tau(w);1)W
                                = \sum_W\widehat{\zeta}(\mathrm{reg}\circ\tau(w))W,                                
\end{align*}
where $\tau=\pi(1-z)\circ S$ is an anti-involution which maps $x \mapsto y, \ y \mapsto x$
and preserves $\mathfrak{h}^0$.
As a consequence, we have
\begin{prop} The Drinfel$'$d associator satisfies \cite{D}
\begin{align}
        \varphi_{\mathrm{KZ}}(X,Y)\cdot\varphi_{\mathrm{KZ}}(-Y,-X) = 1,
\end{align}
which is equivalent to a generalization of the duality formula in \cite{zagier}:
\begin{align}
 \widehat{\zeta}(\mathrm{reg}(w)) =  \widehat{\zeta}(\mathrm{reg}\circ\tau(w)),
\end{align}
where $w $ is an arbitrary word. 
\end{prop} 
(This equivalency was firstly pointed out to the authors by Masanobu 
Kaneko in private communications.)

Moreover writing down \eqref{eq:full euler}, we have
\begin{align*}
 \sum_W \left(\sum_{w_1w_2=w} \mathrm{Li}(\tau(w_1);1-z)\mathrm{Li}(w_2;z)\right) W
 =
 \sum_W \widehat{\zeta}(\mathrm{reg}(w)) W,
\end{align*}
which is equivalent to the Euler connection formula for multiple polylogarithms:
\begin{prop} We have
\begin{align}
 \sum_{w_1w_2=w} \mathrm{Li}(\tau(w_1);1-z)\mathrm{Li}(w_2;z)
 = \widehat{\zeta}(\mathrm{reg}(w)).
\end{align}
where $w$ is an arbitrary word. 
\end{prop}

Let us discuss more on the connection problem of \eqref{FKZ}. The hexagon relations
in \cite{D} can be derived from a braid relation
\begin{align}
   \frac{1}{z}  = (1-z) \circ \left(\frac{z}{z-1}\right) \circ (1-z)
                           = \left(\frac{z}{z-1}\right) \circ (1-z) \circ \left(\frac{z}{z-1}\right).
\end{align}
First applying $\pi(1/z)$ to $G_0$, we have
\begin{align*}
 G_{\infty} &= \pi\left(\frac{1}{z}\right)(G_0)
 =
  \pi(1-z) \circ \pi\left(\frac{z}{z-1}\right) \circ \pi(1-z) (G_0)\\
 &=
  \pi(1-z) \circ \pi\left(\frac{z}{z-1}\right)(G_0 \times \varphi_{\mathrm{KZ}}(X,Y)^{-1})\\
 &=
  \pi(1-z)(G_0\times\exp(\mp X\pi i)) \times \varphi_{\mathrm{KZ}}(X,X-Y)^{-1}\\
 &=
  G_0\times\varphi(X,Y)^{-1}\times\exp(\mp Y\pi i) \times \varphi_{\mathrm{KZ}}(-Y,X-Y)^{-1}.
\end{align*}
By similar computation, we have
\begin{align*}
 G_\infty &= \pi\left(\frac{1}{z}\right)(G_0)
 =
  \pi\left(\frac{z}{z-1}\right) \circ \pi(1-z) \circ 
                               \pi\left(\frac{z}{z-1}\right)(G_0)\\
 &=
  \pi\left(\frac{z}{z-1}\right) \circ \pi(1-z)(G_0\times\exp(\mp X\pi i))\\
 &=
  \pi\left(\frac{z}{z-1}\right)
                  (G_0\times\varphi_{\mathrm{KZ}}(X,Y)^{-1}\times\exp(\mp Y\pi i))\\
 &=
  G_0 \times \exp(\mp X\pi i)
   \times\varphi_{\mathrm{KZ}}(X,X-Y)^{-1}\times\exp(\mp(-X+Y)\pi i).
\end{align*}
Therefore we have the hexagon relations:
\begin{multline}
\exp(\pm X\pi i)
=
 \varphi_{\mathrm{KZ}}(-X+Y,-X)
  \times
 \exp(\mp(-X+Y)\pi i)\\
  \times
 \varphi_{\mathrm{KZ}}(-X+Y,Y)^{-1}
  \times
 \exp(\pm Y\pi i)
  \times 
 \varphi_{\mathrm{KZ}}(X,Y).
\end{multline}

\section{Mellin transforms of polylogarithms and the sum formula for MZVs}

The Mellin transform and the inverse Mellin transform were used in [OU] to show
that the reduced Ohno relation is converted to the Landen connection formula 
for MPLs. In this section, we will consider the Mellin transforms 
of the Euler connection formula and the Landen connection formula for polylogarithms.

We define the Mellin transform by
\begin{align}\label{MELL}
       M[\varphi(z)](\lambda) = \int_0^1\varphi(z)z^{\lambda-1}\,dz,
\end{align}
the inverse Mellin transform by
\begin{align}\label{IMELL}
      \widetilde{M}[f(\lambda)](z) = \frac{1}{2\pi \sqrt{-1}}\int_C
                                     f(\lambda)z^{\lambda}\,d\lambda.
\end{align}
For the details (the integral contour C in (\ref{IMELL}), and the relation 
between both transforms, for example), see [OU].

Now we calculate the Mellin transform of the Euler connection formula 
for the polylogarithms, 
\begin{align}\label{EULPOL}
\mathrm{Li}_k(z) + \mathrm{Li}_{k-1}(z)\mathrm{Li}_{1}(1-z) + \dotsb
  +\mathrm{Li}_{1}(z)\mathrm{Li}_{\underbrace{\scriptstyle1\dotsc1}_{k-1}}(1-z)+ \nonumber \\
  +\mathrm{Li}_{\underbrace{\scriptstyle 21\dotsc1}_{k-1}}(1-z)
 = \zeta(k).
\end{align}
Since $\mathrm{Li}_k(z)=\sum_{n=1}^{\infty}\frac{z^n}{n^k} \ (k\geq 2)$ is uniformly
convergent for $0\leq z \leq 1$, and 
\begin{align*}
\mathrm{Li}_{\underbrace{\scriptstyle 1\dotsc1}_{j}}(1-z) z^{-\lambda-1} =\frac{1}{j!}(-\log z)^j z^{-\lambda-1}
          = \frac{1}{j!}\left(\frac{d}{d\lambda}\right)^jz^{-\lambda-1},
\end{align*}
we have
\begin{align*}
  M[\mathrm{Li}_{k-j}(z) \mathrm{Li}_{\underbrace{\scriptstyle 1\dotsc1}_j}(1-z)](\lambda)&=
 \frac{1}{j!}\left(\frac{d}{d\lambda}\right)^j
 \int_0^1
  \mathrm{Li}_{k-j}(z) z^{-\lambda-1}\,dz\\
 &=
 \frac{1}{j!}\left(\frac{d}{d\lambda}\right)^j
  \sum_{n=1}^{\infty}\frac{1}{n^{k-j}(n-\lambda)} \\
  &=
\sum_{n=1}^{\infty}\frac{1}{n^{k-j}(n-\lambda)^{j+1}}.
\end{align*}
Hence the sum 
$\sum_{j=0}^{k-1}M[\mathrm{Li}_{k-j}(z) \mathrm{Li}_{\underbrace{\scriptstyle 1\dotsc1}_j}(1-z)](\lambda)$
is shown to be
\begin{align*}
   \sum_{j=0}^{k-1}\sum_{n=1}^{\infty} \frac{1}{n^{k-j}(n-\lambda)^{j+1}}
   &=
   \sum_{n=1}^{\infty}\frac{1}{n^k}
         \frac{\dfrac{1}{n-\lambda}-\dfrac{n^k}{(n-\lambda)^{k+1}}}{1-\dfrac{n}{n-\lambda}}\\
   &=
   \sum_{n=1}^{\infty} \frac{1}{-\lambda} \left\{ \frac{1}{n^k} - 
                                          \frac{1}{(n-\lambda)^{k}} \right\} \\
   &=
    -\frac{\zeta(k)}{\lambda}+ \frac{1}{\lambda}
   \sum_{m=0}^{\infty}\binom{-k}{m} \sum_{n=1}^{\infty}\frac{(-\lambda)^m}{n^{k+m}}\\
   &=
   \frac{\zeta(k)}{\lambda}
    + \frac{1}{\lambda} \sum_{m=0}^\infty \binom{m+k-1}{m} \zeta(k+m)\lambda^m.
\end{align*}
The Mellin transform of the last term in \eqref{EULPOL} can be calculated as follows:
\begin{align*}
 & \qquad \int_0^1 \mathrm{Li}_{\underbrace{\scriptstyle 21\dotsc1}_{k-1}}(1-z)z^{-\lambda-1} \,dz
 = \sum_{n_1>\dotsb>n_{k-1}>0} \frac{1}{n_1^2n_2\dotsb n_{k-1}}
  \int_0^1 z^{-\lambda-1} (1-z)^{n_1} \,dz \\
 \intertext{(by the integral representation of the beta function)} 
 &\qquad = \sum_{n_1>\dotsb>n_{k-1}>0} \frac{1}{n_1^2n_2\dotsb n_{k-1}}
  \frac{\Gamma(-\lambda)\Gamma(n_1+1)}{\Gamma(n_1-\lambda+1)}\\
 &\qquad = \frac{1}{-\lambda}\sum_{n_1>\dotsb>n_{k-1}>0} \frac{1}{n_1^2n_2\dotsb n_{k-1}}
       \frac{n_1!}{(n_1-\lambda)(n_1-1-\lambda)\dotsb(1-\lambda)}\\
 &\qquad = \frac{1}{-\lambda}
  \sum_{n_1>\dotsb>n_{k-1}}
    \sum_{m=0}^\infty
    \sum_{\substack{l_i\ge0\\l_1+\dotsb+l_{n_1}=m}}
     \frac{1}{n_1^2n_2\dotsb n_{k-1}}
       \frac{\lambda^m}{{n_1}^{l_{n_1}} (n_1-1)^{l_{n_1-1}}\dotsm 2^{l_2}1^{l_1}}
    \\
 \intertext{(by some combinatorial consideration)}
 &\qquad = \frac{1}{-\lambda} \sum_{m=0}^\infty \sum_{j=0}^{m}
                   \binom{(k-1)+j-1}{j} S(k+m,(k-1)+j) \lambda^m\\
 &\qquad = \frac{1}{-\lambda} \sum_{m=0}^\infty \sum_{j=0}^{m}
                             \binom{k+j-2}{j} S(k+m,(k-1)+j)\lambda^m,
\end{align*}
where $S(n,r)$ denotes the sum of all MZVs of weight $n$ and depth $r$. Thus we obtain,
for $k\ge2$, $m\ge0$,
\begin{align}
  \sum_{j=0}^{m} \binom{k+j-2}{k-2} S(k+m,k-1+j)
 = 
  \binom{m+k-1}{k-1} \zeta(k+m).
\end{align}
which is rewritten, by putting $n=k+m$, $d=k-1+j$, as
\begin{align}\label{sum1}
\sum_{d=k-1}^{n-1} \binom{d-1}{k-2} S(n,d)
 = 
  \binom{n-1}{k-1} \zeta(n)  \quad (2 \leq k \leq n). \tag{$\heartsuit$}
\end{align}
From this we can easily prove the sum formula
\begin{align}\label{SUMF}
  S(n,r) = \zeta(n) \quad (n \geq 2, \  n \geq r \geq 1)
\end{align}
by induction on $r$.
\begin{proof}
That $S(n,n-1) = \zeta(n)$ follows from (\ref{sum1}) with $k=n$. 
For $r<n-1$, we write (\ref{sum1}) as
\begin{align*}
   S(n,r) +
  \sum_{d=r+1}^{n-1} \binom{d-1}{r-1} S(n,d)
 = 
  \binom{n-1}{r} \zeta(n).
\end{align*}
From the induction hypothesis, we have
\begin{align*}
   S(n,r) +
  \sum_{d=r+1}^{n-1} \binom{d-1}{r-1} \zeta(n)
 = 
  \binom{n-1}{r} \zeta(n)
\end{align*}
Noting that
\begin{align*}
  \sum_{d=r+1}^{n-1} \binom{d-1}{r-1}
 = \binom{n-1}{r} - 1,
\end{align*}
we obtain (\ref{SUMF}).
\end{proof}

Next we consider the Mellin transform of the Landen connection formula for polylogarithms:
\begin{align}\label{LANDEN}
 \mathrm{Li}_m(z) = -\sum_{j=1}^m \sum_{\substack{\mathbf{c}\\ \mathrm{weight}\,m\\ \mathrm{length}\,m-j+1}}
     \mathrm{Li}_{\mathbf{c}}\left( \frac{z}{z-1} \right).
\end{align}
For $m=2$, it is nothing but the Landen formula for the dilogarithm 
\begin{align}\label{LANDEND}
    \mathrm{Li}_2(z) = -\mathrm{Li}_2\left( \frac{z}{z-1} \right) - \mathrm{Li}_{11}\left( \frac{z}{z-1} \right).
\end{align}
First we compute the Mellin transform of this formula. The Taylor expansion of (\ref{LANDEND})
of both sides reads
\begin{align*}
  \sum_{n=1}^{\infty} \frac{z^n}{n^2} = \sum_{n_1>n_2>0}\frac{z^{n_2}}{n_1(n_1-n_2)}+
                               \sum_{n_1>n_2>0}\frac{z^{n_1}}{n_1(n_2-n_1)}.
\end{align*}
Applying the Mellin transform to both sides above, we have
\begin{align*}
\sum_{n=1}^{\infty}\frac{1}{n^2(n-\lambda)}&=\sum_{n_1>n_2>0}\frac{1}{n_1(n_1-n_2)(n_2-\lambda)}
                                     + \sum_{n_1>n_2>0}\frac{1}{n_1(n_2-n_1)(n_1-\lambda)}\\
                       &= \sum_{n_1>n_2>0}\frac{1}{n_1(n_1-\lambda)(n_2-\lambda)}.
\end{align*}
The Taylor expansion of the both sides above at $\lambda=0$ reads
\begin{align*}
     \sum_{l=0}^{\infty}\zeta(3+l)\lambda^l=\sum_{l=0}^{\infty} \left\{
           \sum_{\substack{c_1+c_2=l\\c_1,c_2\geq0}}\zeta(2+c_1,1+c_2) \right\}\lambda^l,
\end{align*}
which yields the sum formula for MZVs of the depth 2;
\begin{align}
       \zeta(3+l) = \sum_{\substack{c_1+c_2=l\\c_1,c_2\geq0}}\zeta(2+c_1,1+c_2).
\end{align}

Now we consider the general case. What we have to show is 
\begin{lemma}[cf. \cite{OU}]  For $1 \leq j \leq m$, 
\begin{align}\label{LANDENLEMMA}
\sum_{n_1>\ldots>n_m>0} \frac{z^{n_j}}{n_1\prod_{i\neq j}(n_i-n_j)} = -
            \sum_{\substack{\mathbf{c}\\ \mathrm{weight}\,m\\ \mathrm{length}\,m-j+1}}
            \mathrm{Li}_{\mathbf{c}}\left( \frac{z}{z-1} \right).
\end{align}
\end{lemma}

\begin{proof}
We prove this by induction on $m$. For $m=1$, it is obvious. Suppose that it holds
for $m-1$. Let $j \neq 1,m$. Then we have
\begin{align*}
&
\frac{d}{dz} \left( \mbox{the LHS of (\ref{LANDENLEMMA})} \right)\\
&=
 \sum_{n_1>\dotsb>n_m>0}
  \frac{n_j z^{n_j-1}}
       {n_1(n_1-n_j)\dotsb(n_{j-1}-n_j)\cdot(n_{j+1}-n_j)\dotsb(n_m-n_j)}
      \tag{$\clubsuit$}
 \end{align*}
 Putting $n_i=l_i+\dotsb+l_m \quad (1 \leq i \leq m)$,
 we have
 \begin{align*}
 &  (\clubsuit) =
 \sum_{l_1,\dotsc,l_m=1}^\infty
  \left\{\frac{1}{l_1+\dotsb+l_{j-1}}-\frac{1}{l_1+\dotsb+l_m}\right\} \times \\
 & \hspace{3cm} \times 
  \frac{(-1)^{m-j}z^{l_j+\dotsb+l_m-1}}
      {(l_2+\dotsb+l_{j-1})\dotsm l_{j-1}\cdot l_j\dotsm(l_j+\dotsb+l_{m-1})}\\
 &=
 \sum_{l_2,\dotsc,l_m=1}^\infty
 \left\{\sum_{l_1=1}^{l_j+\dotsb+l_m}\frac{1}{l_1+\dotsb+l_{j-1}}\right\}
  \frac{(-1)^{m-j}z^{l_j+\dotsb+l_m-1}}
      {(l_2+\dotsb+l_{j-1})\dotsm l_{j-1}\cdot l_j\dotsm(l_j+\dotsb+l_{m-1})}\\
 &=
 \sum_{l_2,\dotsc,l_{m-1}=1}^\infty
 \left\{
  \sum_{l_1=1}^{l_j+\dotsb+l_{m-1}} \sum_{l_m=1}^\infty
  + \sum_{l_1=l_j+\dotsb+l_{m-1}+1}^\infty
    \cdot \sum_{l_m=l_1-(l_j+\dotsb+l_{m-1})}^\infty\right\}\\
 &  \hspace{3cm} \frac{(-1)^{m-j}z^{l_j+\dotsb+l_m-1}}
   {(l_1+\cdots+l_{j-1})(l_2+\dotsb+l_{j-1})\dotsm l_{j-1}\cdot l_j\dotsm(l_j+\dotsb+l_{m-1})}\\
 &=
 \frac{1}{1-z}
 \sum_{l_1,\dotsc,l_{m-1}=1}^\infty\sum_{l_1=1}^{l_j+\dotsb+l_{m-1}}
  \frac{(-1)^{m-j}z^{l_j+\dotsb+l_{m-1}}}{\mbox{\{the same denominator as above\}}}\\
      %{(l_1+\cdots+l_{j-1})(l_2+\cdots+l_{j-1})\cdots l_{j-1}
       % \cdot l_j\cdots(l_j+\cdots+l_{m-1})}\\
 &\quad
 + \frac{1}{1-z}
 \sum_{l_1,\dotsc,l_{m-1}=1}^\infty \sum_{l_1=l_j+\dotsb+l_{m-1}+1}^\infty
  \frac{(-1)^{m-j}z^{l_1-1}}{\mbox{\{the same denominator as above\}}}\\
      %{(l_1+\cdots+l_{j-1})(l_2+\cdots+l_{j-1})\cdots l_{j-1}
        %\cdot l_j\dotsb(l_j+\cdots+l_{m-1})}  \\
 &=
 \frac{1}{1-z}
 \sum_{l_1,\dotsc,l_{m-1}=1}^\infty
  \left\{
   \frac{1}{l_1+\dotsb+l_{j-1}} - \frac{1}{l_1+\dotsb+l_{m-1}}
  \right\} \times \\
 & \hspace{6cm}
  \frac{(-1)^{m-j}z^{l_j+\dotsb+l_{m-1}}}
      {(l_2+\dotsb+l_{j-1})\dotsm l_{j-1}
        \cdot l_j\dotsm(l_j+\dotsb+l_{m-1})}\\
 &\quad
 + \frac{1}{1-z}
  \sum_{l_1,\dotsc,l_{m-1}=1}^\infty
  \frac{(-1)^{m-j}z^{l_j+\dotsb+l_m-1}}
      {(l_2+\dotsb+l_m)(l_2+\dotsb+l_{j-1})\dotsm l_{j-1}
        \cdot l_j\dotsm(l_j+\dotsb+l_{m-1})}\\
 &=
 \frac{1}{1-z}
 \sum_{l_1,\dotsc,l_{m-1}=1}^\infty
  \frac{(-1)^{m-j}z^{l_j+\dotsb+l_{m-1}}}
      {(l_1+\dotsb+l_{m-1})(l_1+\dotsb+l_{j-1})\dotsm l_{j-1}
        \cdot l_j\dotsm(l_j+\dotsb+l_{m-2})}\\
 &
 + \frac{1}{z(1-z)}
  \sum_{l_1,\dotsc,l_{m-1}=1}^\infty
  \frac{(-1)^{m-j}z^{l_{j-1}+\dotsb+l_{m-1}}}
      {(l_1+\dotsb+l_{m-1})(l_1+\dotsb+l_{j-2})\dotsm l_{j-2}
        \cdot l_{j-1}\dotsm(l_{j-1}+\dotsb+l_{m-2})}\\
 \intertext{(putting $n_i=l_i+\cdots+l_{m-1} \quad (1 \leq i \leq m-1)$)}
 & \quad
 = - \frac{1}{1-z}\sum_{n_1>\dotsb>n_{m-1}>0}
  \frac{z^{n_j}}{n_1(n_1-n_j)\dotsm(n_{j-1}-n_{j})
   \cdot (n_{j+1}-n_{j})\dotsm(n_{m-1}-n_{j})}\\
 &\quad
 + \frac{1}{z(1-z)}\sum_{n_1>\dotsb>n_{m-1}>0}
    \frac{z^{n_{j-1}}}{n_1(n_1-n_{j-1})\dotsm(n_{j-2}-n_{j-1})
     \cdot (n_{j}-n_{j-1})\dotsm(n_{m-1}-n_{j-1})}\\
\intertext{(by the induction hypothesis)}
& \quad
= \frac{1}{1-z} \sum_{\substack{\mathbf{c}\\ \mathrm{weight}\, m-1\\ \mathrm{length}\,m-j}}
            \mathrm{Li}_{\mathbf{c}}\left(\frac{z}{z-1}\right) 
-\frac{1}{z(1-z)}
               \sum_{\substack{\mathbf{c}\\ \mathrm{weight}\, m-1\\ \mathrm{length}\,m-j+1}}
            \mathrm{Li}_{\mathbf{c}}\left(\frac{z}{z-1}\right).
\end{align*}
Using the differential relation,
\begin{align*}
\frac{d}{dz}\mathrm{Li}_{k_1,\ldots,k_m}\left(\frac{z}{z-1}\right)=
\begin{cases}
 \dfrac{1}{z(1-z)}\mathrm{Li}_{k_1-1,k_2,\ldots,k_m}\left(\dfrac{z}{z-1}\right) & \mbox{if} \quad k_1\geq 2, \\
{}  & \\
 - \dfrac{1}{1-z}\mathrm{Li}_{k_2,\ldots,k_m}\left(\dfrac{z}{z-1}\right) & \mbox{if} \quad k_1=1
\end{cases}
\end{align*}
we obtain (\ref{LANDENLEMMA}) for $j \neq 1,m$. 
The case for $j=1,m$ can be shown in the same way.
\end{proof}
From the Landen formula (\ref{LANDEN}), we have
\begin{align*}
     \sum_{n=1}^{\infty}\frac{z^n}{n^m} =
     \sum_{j=1}^m \sum_{n_1>\dotsb>n_m>0} \frac{z^{n_j}}{n_1\prod_{i\neq j}(n_i-n_j)}.
\end{align*}
The Mellin transform of the above reads 
\begin{align*}
     \sum_{n=1}^{\infty}\frac{1}{n^m(n-\lambda)} =
     \sum_{j=1}^m \sum_{n_1>\dotsb>n_m>0} \frac{1}{n_1\prod_{i\neq j}(n_i-n_j)(n_j-\lambda)}.
\end{align*}
The Taylor expansion of the both sides yields the sum formula (\ref{SUMF}).

\section*{Appendix 1: 
Knizhnik-Zamolodchikov equation over the configuration space $\mathbb{F}_3(\mathbb{C})$}

Let $\mathfrak{P}_n
=\mathbb{C}\langle\langle X_{ij}\rangle\rangle_{1\le i \neq j \le n}$ be an algebra
with the defining relations
\begin{align}\label{DEFRELATION}
 \begin{cases}
 X_{ij} = X_{ji},& \\
 [ X_{ij}, X_{kl}] =0 ,     & \text{($i,j,k,l$ distinct)}\\
 [ X_{ij}, X_{ik} + X_{kj}]=0. & \text{($i,j,k$ distinct)}
 \end{cases}
\end{align}

Putting $z_{ij} := z_i-z_j$, we consider the following system of differential
equations:
\begin{align}
 dW = \left( \sum_{i<j} X_{ij} d\log z_{ij} \right) W. \tag{KZ$_n$}
\end{align}
This system is integrable because of \eqref{DEFRELATION} (cf.\ \cite{Kas}).

In what follows, we consider only the case of $n=3$.
We set $T = X_{12} + X_{13} + X_{23}$ which belongs to the center of $\mathfrak{P}_3$.
\begin{prop}
Solutions to \eqref{KZ3} satisfying the following asymptotic behavior exist, 
\begin{align}\label{ASYMPTOTIC}
 W_{(ij)k} &\sim z_{ji}^{X_{ij}} z_{ki}^{X_{ik}+X_{jk}}
   \qquad (|z_{ji}| << |z_{ki}|),\\
 W_{i(jk)} &\sim z_{kj}^{X_{jk}} z_{ki}^{X_{ij}+X_{ik}}
   \qquad (|z_{kj}| << |z_{ki}|),
\end{align}
and are expressed in terms of the solutions $G_0$ and $G_1$ of \eqref{FKZ} as follows$:$
\begin{align}\label{KZ3SOL}
 W_{(ij)k}(z_1, z_2, z_3)
 &= G_0\left(X_{ij},-X_{jk};\frac{z_{ji}}{z_{ki}}\right)
    \times z_{ki}^{T},\\
W_{i(jk)}(z_1, z_2, z_3)
 &= G_1\left(X_{ij},-X_{jk};\frac{z_{ji}}{z_{ki}}\right)
    \times z_{ki}^T.
\end{align}
\end{prop}

\begin{proof}
We set
\begin{align*}
 W(z_1, z_2, z_3)
= G\left(\frac{z_{ji}}{z_{ki}}\right) z_{ki}^T.
\end{align*}
Then $G$ satisfies 
\begin{align*}
 dG\left(\frac{z_{ji}}{z_{ki}}\right) z_{ki}^T + 
 G\left(\frac{z_{ji}}{z_{ki}}\right) z_{ki}^T T d\log z_{ki}
 = 
 \left( \sum_{p<q} X_{pq} d\log z_{pq} \right)
 G\left(\frac{z_{ji}}{z_{ki}}\right) z_{ki}^T,
\end{align*}
so that we have
\begin{align}
 dG \left(\frac{z_{ji}}{z_{ki}}\right)
 = \left(
    X_{ij}d\log \frac{z_{ji}}{z_{ki}}
    + X_{jk}d\log\left( 1 - \frac{z_{ji}}{z_{ki}} \right)
   \right) G \left(\frac{z_{ji}}{z_{ki}}\right).
\end{align}
Taking into account \eqref{ASYMPTOTIC}, we obtain the formula \eqref{KZ3SOL}.
\end{proof}

\begin{prop}
The solutions above satisfy the connection formulas
\begin{align}
& \quad   W_{i(jk)}(z_1, z_2, z_3) \times \varphi_{\mathrm{KZ}}(X_{ij},-X_{jk})
             = W_{(ij)k}(z_1, z_2, z_3),  \label{CON}\\
& \quad  W_{(ji)k}(z_1, z_2, z_3) \times \exp(\pm X_{ij} \pi i)
    = W_{(ij)k}(z_1,z_2,z_3) \qquad \frac{z_{ji}}{z_{ki}} \in \mathbf{H}_{\pm}, \label{CONNE} \\
& \quad W_{i(kj)}(z_1, z_2, z_3) \times \exp(\pm X_{jk} \pi i)
    = W_{i(jk)}(z_1, z_2, z_3) \qquad \frac{z_{ji}}{z_{ki}} \in \mathbf{H}_{\mp}, \label{CONNEC}\\
& \quad W_{k(ij)}(z_1,z_2,z_3)  \times \exp(\pm (X_{ik}+X_{jk})\pi i)
    = W_{(ij)k}(z_1, z_2, z_3) \qquad \frac{z_{ji}}{z_{ki}} \in \mathbf{H}_{\mp}. \label{CONNECT}\\
& \quad W_{(jk)i}(z_1,z_2,z_3)  \times \exp(\pm (X_{ik}+X_{jk})\pi i)
    = W_{i(jk)}(z_1, z_2, z_3) \qquad \frac{z_{ji}}{z_{ki}} \in \mathbf{H}_{\pm}. 
 \label{CONNECTI}
\end{align}
\end{prop}

Here in each equation,
$\frac{z_{ji}}{z_{ki}}\in\mathbf{H}_{\pm}$ corresponds to
the analytic continuation of the paths Figure~1-4:

\vspace{.5cm}
\begin{figure}[p]
\unitlength 1mm
\begin{picture}(105.00,15.63)(0,0)

\thicklines
\path(10.00,5.63)(105.00,5.63)

\put(94.38,5.63){\makebox(0,0)[cc]{$\bullet$}}

\put(55.00,5.63){\makebox(0,0)[cc]{$\bullet$}}

\put(44.38,5.63){\makebox(0,0)[cc]{$\bullet$}}

\put(35.00,5.63){\makebox(0,0)[cc]{$\bullet$}}

\thicklines
\put(45.00,5.63){\arc{20.00}{-2.97}{-0.17}}

\thicklines
\path(35.00,7.50)(34.38,8.75)

\thicklines
\path(35.00,7.50)(36.25,8.75)

\put(70.63,15.63){\makebox(0,0)[cc]{$\frac{z_{ji}}{z_{ki}}\in\mathbf{H}_{+}$}}

\put(60.63,15.00){\makebox(0,0)[cc]{}}

\put(35.00,1.88){\makebox(0,0)[cc]{$z_j$}}

\put(55.00,1.88){\makebox(0,0)[cc]{$z_{j}$}}

\put(45.00,1.88){\makebox(0,0)[cc]{$z_{i}$}}

\put(95.00,1.88){\makebox(0,0)[cc]{$z_{k}$}}

\end{picture}
 \caption{$z_j$ goes around $z_i$ by counter-clockwise}
 \label{fig:ij}
\end{figure}

\begin{figure}[p]
\unitlength 1mm
\begin{picture}(105.00,15.63)(0,0)

\thicklines
\path(10.00,5.63)(105.00,5.63)

\put(95.63,5.63){\makebox(0,0)[cc]{$\bullet$}}

\put(85.00,5.63){\makebox(0,0)[cc]{$\bullet$}}

\put(76.25,5.63){\makebox(0,0)[cc]{$\bullet$}}

\put(35.00,5.63){\makebox(0,0)[cc]{$\bullet$}}

\thicklines
\put(85.63,5.63){\arc{20.00}{-2.97}{-0.17}}

\thicklines
\path(75.63,7.50)(75.00,8.75)

\thicklines
\path(75.63,7.50)(76.88,8.75)

\put(51.25,15.63){\makebox(0,0)[cc]{$\frac{z_{ji}}{z_{ki}}\in\mathbf{H}_{-}$}}

\put(60.63,15.00){\makebox(0,0)[cc]{}}

\put(35.00,1.88){\makebox(0,0)[cc]{$z_i$}}

\put(85.63,2.50){\makebox(0,0)[cc]{$z_{j}$}}

\put(75.63,2.50){\makebox(0,0)[cc]{$z_{k}$}}

\put(95.00,1.88){\makebox(0,0)[cc]{$z_{k}$}}

\end{picture}

 \caption{$z_k$ goes around $z_j$ by counter-clockwise}
 \label{fig:jk}
\end{figure}

\begin{figure}[p]
\unitlength 1mm
\begin{picture}(105.00,16.88)(0,0)

\thicklines
\path(10.00,5.63)(105.00,5.63)

\put(95.00,5.63){\makebox(0,0)[cc]{$\bullet$}}

\put(65.00,5.63){\makebox(0,0)[cc]{$\bullet$}}

\put(31.88,5.63){\makebox(0,0)[cc]{$\bullet$}}

\put(55.00,5.63){\makebox(0,0)[cc]{$\bullet$}}

\thicklines
\path(32.50,6.88)(31.88,8.13)

\thicklines
\path(32.50,6.88)(33.75,8.13)

\put(21.25,15.00){\makebox(0,0)[cc]{$\frac{z_{ji}}{z_{ki}}\in\mathbf{H}_{-}$}}

\put(60.63,15.00){\makebox(0,0)[cc]{}}

\put(54.38,2.50){\makebox(0,0)[cc]{$z_i$}}

\put(65.00,2.50){\makebox(0,0)[cc]{$z_{j}$}}

\put(31.25,2.50){\makebox(0,0)[cc]{$z_{k}$}}

\put(95.00,1.88){\makebox(0,0)[cc]{$z_{k}$}}

\thicklines
\multiput(94.78,7.50)(0.10,-0.34){1}{\line(0,-1){0.34}}
\multiput(94.65,7.84)(0.13,-0.34){1}{\line(0,-1){0.34}}
\multiput(94.49,8.17)(0.16,-0.33){1}{\line(0,-1){0.33}}
\multiput(94.29,8.50)(0.10,-0.17){2}{\line(0,-1){0.17}}
\multiput(94.06,8.83)(0.11,-0.17){2}{\line(0,-1){0.17}}
\multiput(93.80,9.16)(0.13,-0.16){2}{\line(0,-1){0.16}}
\multiput(93.52,9.49)(0.14,-0.16){2}{\line(0,-1){0.16}}
\multiput(93.20,9.81)(0.11,-0.11){3}{\line(0,-1){0.11}}
\multiput(92.85,10.13)(0.12,-0.11){3}{\line(1,0){0.12}}
\multiput(92.47,10.44)(0.13,-0.10){3}{\line(1,0){0.13}}
\multiput(92.06,10.75)(0.14,-0.10){3}{\line(1,0){0.14}}
\multiput(91.62,11.05)(0.15,-0.10){3}{\line(1,0){0.15}}
\multiput(91.16,11.35)(0.16,-0.10){3}{\line(1,0){0.16}}
\multiput(90.67,11.65)(0.25,-0.15){2}{\line(1,0){0.25}}
\multiput(90.15,11.94)(0.26,-0.14){2}{\line(1,0){0.26}}
\multiput(89.60,12.22)(0.27,-0.14){2}{\line(1,0){0.27}}
\multiput(89.02,12.50)(0.29,-0.14){2}{\line(1,0){0.29}}
\multiput(88.42,12.77)(0.30,-0.14){2}{\line(1,0){0.30}}
\multiput(87.80,13.03)(0.31,-0.13){2}{\line(1,0){0.31}}
\multiput(87.15,13.29)(0.32,-0.13){2}{\line(1,0){0.32}}
\multiput(86.48,13.54)(0.34,-0.13){2}{\line(1,0){0.34}}
\multiput(85.78,13.79)(0.35,-0.12){2}{\line(1,0){0.35}}
\multiput(85.06,14.02)(0.36,-0.12){2}{\line(1,0){0.36}}
\multiput(84.32,14.25)(0.37,-0.11){2}{\line(1,0){0.37}}
\multiput(83.56,14.47)(0.38,-0.11){2}{\line(1,0){0.38}}
\multiput(82.78,14.68)(0.39,-0.11){2}{\line(1,0){0.39}}
\multiput(81.98,14.88)(0.40,-0.10){2}{\line(1,0){0.40}}
\multiput(81.16,15.07)(0.41,-0.10){2}{\line(1,0){0.41}}
\multiput(80.32,15.26)(0.42,-0.09){2}{\line(1,0){0.42}}
\multiput(79.47,15.43)(0.85,-0.18){1}{\line(1,0){0.85}}
\multiput(78.60,15.60)(0.87,-0.17){1}{\line(1,0){0.87}}
\multiput(77.71,15.76)(0.89,-0.16){1}{\line(1,0){0.89}}
\multiput(76.81,15.90)(0.90,-0.15){1}{\line(1,0){0.90}}
\multiput(75.90,16.04)(0.91,-0.14){1}{\line(1,0){0.91}}
\multiput(74.97,16.17)(0.93,-0.13){1}{\line(1,0){0.93}}
\multiput(74.04,16.28)(0.94,-0.12){1}{\line(1,0){0.94}}
\multiput(73.09,16.39)(0.95,-0.11){1}{\line(1,0){0.95}}
\multiput(72.13,16.49)(0.96,-0.10){1}{\line(1,0){0.96}}
\multiput(71.17,16.57)(0.97,-0.09){1}{\line(1,0){0.97}}
\multiput(70.20,16.65)(0.97,-0.08){1}{\line(1,0){0.97}}
\multiput(69.22,16.71)(0.98,-0.06){1}{\line(1,0){0.98}}
\multiput(68.23,16.77)(0.98,-0.05){1}{\line(1,0){0.98}}
\multiput(67.24,16.81)(0.99,-0.04){1}{\line(1,0){0.99}}
\multiput(66.25,16.84)(0.99,-0.03){1}{\line(1,0){0.99}}
\multiput(65.25,16.86)(1.00,-0.02){1}{\line(1,0){1.00}}
\multiput(64.26,16.87)(1.00,-0.01){1}{\line(1,0){1.00}}
\put(63.26,16.87){\line(1,0){1.00}}
\multiput(62.26,16.86)(1.00,0.01){1}{\line(1,0){1.00}}
\multiput(61.26,16.84)(1.00,0.02){1}{\line(1,0){1.00}}
\multiput(60.27,16.81)(0.99,0.03){1}{\line(1,0){0.99}}
\multiput(59.28,16.77)(0.99,0.04){1}{\line(1,0){0.99}}
\multiput(58.30,16.71)(0.99,0.05){1}{\line(1,0){0.99}}
\multiput(57.32,16.65)(0.98,0.06){1}{\line(1,0){0.98}}
\multiput(56.34,16.57)(0.97,0.08){1}{\line(1,0){0.97}}
\multiput(55.38,16.49)(0.97,0.09){1}{\line(1,0){0.97}}
\multiput(54.42,16.39)(0.96,0.10){1}{\line(1,0){0.96}}
\multiput(53.48,16.28)(0.95,0.11){1}{\line(1,0){0.95}}
\multiput(52.54,16.17)(0.94,0.12){1}{\line(1,0){0.94}}
\multiput(51.61,16.04)(0.93,0.13){1}{\line(1,0){0.93}}
\multiput(50.70,15.90)(0.91,0.14){1}{\line(1,0){0.91}}
\multiput(49.80,15.76)(0.90,0.15){1}{\line(1,0){0.90}}
\multiput(48.91,15.60)(0.89,0.16){1}{\line(1,0){0.89}}
\multiput(48.04,15.44)(0.87,0.17){1}{\line(1,0){0.87}}
\multiput(47.19,15.26)(0.85,0.18){1}{\line(1,0){0.85}}
\multiput(46.35,15.08)(0.42,0.09){2}{\line(1,0){0.42}}
\multiput(45.53,14.88)(0.41,0.10){2}{\line(1,0){0.41}}
\multiput(44.73,14.68)(0.40,0.10){2}{\line(1,0){0.40}}
\multiput(43.95,14.47)(0.39,0.11){2}{\line(1,0){0.39}}
\multiput(43.19,14.25)(0.38,0.11){2}{\line(1,0){0.38}}
\multiput(42.45,14.02)(0.37,0.11){2}{\line(1,0){0.37}}
\multiput(41.73,13.79)(0.36,0.12){2}{\line(1,0){0.36}}
\multiput(41.03,13.55)(0.35,0.12){2}{\line(1,0){0.35}}
\multiput(40.36,13.30)(0.34,0.13){2}{\line(1,0){0.34}}
\multiput(39.71,13.04)(0.32,0.13){2}{\line(1,0){0.32}}
\multiput(39.08,12.77)(0.31,0.13){2}{\line(1,0){0.31}}
\multiput(38.48,12.50)(0.30,0.14){2}{\line(1,0){0.30}}
\multiput(37.91,12.23)(0.29,0.14){2}{\line(1,0){0.29}}
\multiput(37.36,11.94)(0.27,0.14){2}{\line(1,0){0.27}}
\multiput(36.84,11.65)(0.26,0.14){2}{\line(1,0){0.26}}
\multiput(36.35,11.36)(0.25,0.15){2}{\line(1,0){0.25}}
\multiput(35.88,11.06)(0.16,0.10){3}{\line(1,0){0.16}}
\multiput(35.44,10.75)(0.15,0.10){3}{\line(1,0){0.15}}
\multiput(35.04,10.44)(0.14,0.10){3}{\line(1,0){0.14}}
\multiput(34.66,10.13)(0.13,0.10){3}{\line(1,0){0.13}}
\multiput(34.31,9.81)(0.12,0.11){3}{\line(1,0){0.12}}
\multiput(33.99,9.49)(0.11,0.11){3}{\line(0,1){0.11}}
\multiput(33.70,9.17)(0.14,0.16){2}{\line(0,1){0.16}}
\multiput(33.44,8.84)(0.13,0.16){2}{\line(0,1){0.16}}
\multiput(33.21,8.51)(0.11,0.17){2}{\line(0,1){0.17}}
\multiput(33.02,8.18)(0.10,0.17){2}{\line(0,1){0.17}}
\multiput(32.85,7.84)(0.17,0.33){1}{\line(0,1){0.33}}
\multiput(32.72,7.50)(0.13,0.34){1}{\line(0,1){0.34}}
\multiput(32.62,7.17)(0.10,0.34){1}{\line(0,1){0.34}}
\multiput(32.55,6.83)(0.07,0.34){1}{\line(0,1){0.34}}

\end{picture}
 \caption{$z_k$ goes around $z_i$ and $z_j$ by counter-clockwise}
 \label{fig:(ij)k}
\end{figure}

\begin{figure}[p]
\unitlength 1mm
\begin{picture}(105.00,19.38)(0,0)

\thicklines
\path(10.00,10.00)(105.00,10.00)

\put(95.00,10.00){\makebox(0,0)[cc]{$\bullet$}}

\put(65.00,10.00){\makebox(0,0)[cc]{$\bullet$}}

\put(30.63,10.00){\makebox(0,0)[cc]{$\bullet$}}

\put(55.00,10.00){\makebox(0,0)[cc]{$\bullet$}}

\put(21.25,19.38){\makebox(0,0)[cc]{$\frac{z_{ji}}{z_{ki}}\in\mathbf{H}_{+}$}}

\put(60.63,19.38){\makebox(0,0)[cc]{}}

\put(53.75,13.13){\makebox(0,0)[cc]{$z_j$}}

\put(65.00,13.13){\makebox(0,0)[cc]{$z_{k}$}}

\put(30.00,13.13){\makebox(0,0)[cc]{$z_{i}$}}

\put(94.38,13.13){\makebox(0,0)[cc]{$z_{i}$}}

\thicklines
\multiput(30.73,8.06)(0.09,-0.27){1}{\line(0,-1){0.27}}
\multiput(30.82,7.79)(0.12,-0.27){1}{\line(0,-1){0.27}}
\multiput(30.94,7.52)(0.15,-0.27){1}{\line(0,-1){0.27}}
\multiput(31.10,7.25)(0.09,-0.13){2}{\line(0,-1){0.13}}
\multiput(31.28,6.99)(0.11,-0.13){2}{\line(0,-1){0.13}}
\multiput(31.50,6.72)(0.12,-0.13){2}{\line(0,-1){0.13}}
\multiput(31.75,6.46)(0.14,-0.13){2}{\line(1,0){0.14}}
\multiput(32.02,6.20)(0.15,-0.13){2}{\line(1,0){0.15}}
\multiput(32.32,5.94)(0.17,-0.13){2}{\line(1,0){0.17}}
\multiput(32.66,5.69)(0.18,-0.13){2}{\line(1,0){0.18}}
\multiput(33.02,5.44)(0.20,-0.12){2}{\line(1,0){0.20}}
\multiput(33.41,5.19)(0.21,-0.12){2}{\line(1,0){0.21}}
\multiput(33.83,4.94)(0.22,-0.12){2}{\line(1,0){0.22}}
\multiput(34.28,4.70)(0.24,-0.12){2}{\line(1,0){0.24}}
\multiput(34.75,4.46)(0.25,-0.12){2}{\line(1,0){0.25}}
\multiput(35.25,4.23)(0.26,-0.11){2}{\line(1,0){0.26}}
\multiput(35.78,4.00)(0.28,-0.11){2}{\line(1,0){0.28}}
\multiput(36.33,3.78)(0.29,-0.11){2}{\line(1,0){0.29}}
\multiput(36.91,3.56)(0.30,-0.11){2}{\line(1,0){0.30}}
\multiput(37.51,3.34)(0.31,-0.11){2}{\line(1,0){0.31}}
\multiput(38.14,3.13)(0.33,-0.10){2}{\line(1,0){0.33}}
\multiput(38.79,2.93)(0.34,-0.10){2}{\line(1,0){0.34}}
\multiput(39.47,2.73)(0.35,-0.10){2}{\line(1,0){0.35}}
\multiput(40.16,2.53)(0.36,-0.09){2}{\line(1,0){0.36}}
\multiput(40.88,2.35)(0.37,-0.09){2}{\line(1,0){0.37}}
\multiput(41.62,2.16)(0.76,-0.17){1}{\line(1,0){0.76}}
\multiput(42.38,1.99)(0.78,-0.17){1}{\line(1,0){0.78}}
\multiput(43.16,1.82)(0.80,-0.16){1}{\line(1,0){0.80}}
\multiput(43.96,1.66)(0.82,-0.16){1}{\line(1,0){0.82}}
\multiput(44.77,1.50)(0.83,-0.15){1}{\line(1,0){0.83}}
\multiput(45.60,1.36)(0.85,-0.14){1}{\line(1,0){0.85}}
\multiput(46.45,1.21)(0.86,-0.13){1}{\line(1,0){0.86}}
\multiput(47.32,1.08)(0.88,-0.13){1}{\line(1,0){0.88}}
\multiput(48.20,0.95)(0.89,-0.12){1}{\line(1,0){0.89}}
\multiput(49.09,0.83)(0.91,-0.11){1}{\line(1,0){0.91}}
\multiput(50.00,0.72)(0.92,-0.10){1}{\line(1,0){0.92}}
\multiput(50.92,0.62)(0.93,-0.10){1}{\line(1,0){0.93}}
\multiput(51.85,0.52)(0.94,-0.09){1}{\line(1,0){0.94}}
\multiput(52.79,0.44)(0.95,-0.08){1}{\line(1,0){0.95}}
\multiput(53.74,0.35)(0.96,-0.07){1}{\line(1,0){0.96}}
\multiput(54.70,0.28)(0.97,-0.06){1}{\line(1,0){0.97}}
\multiput(55.66,0.22)(0.97,-0.06){1}{\line(1,0){0.97}}
\multiput(56.64,0.16)(0.98,-0.05){1}{\line(1,0){0.98}}
\multiput(57.62,0.11)(0.98,-0.04){1}{\line(1,0){0.98}}
\multiput(58.60,0.08)(0.99,-0.03){1}{\line(1,0){0.99}}
\multiput(59.59,0.04)(0.99,-0.02){1}{\line(1,0){0.99}}
\multiput(60.58,0.02)(0.99,-0.01){1}{\line(1,0){0.99}}
\multiput(61.58,0.01)(1.00,-0.01){1}{\line(1,0){1.00}}
\multiput(62.57,0.00)(1.00,0.00){1}{\line(1,0){1.00}}
\multiput(63.57,0.00)(0.99,0.01){1}{\line(1,0){0.99}}
\multiput(64.56,0.01)(0.99,0.02){1}{\line(1,0){0.99}}
\multiput(65.55,0.03)(0.99,0.03){1}{\line(1,0){0.99}}
\multiput(66.54,0.06)(0.99,0.04){1}{\line(1,0){0.99}}
\multiput(67.53,0.09)(0.98,0.04){1}{\line(1,0){0.98}}
\multiput(68.51,0.14)(0.98,0.05){1}{\line(1,0){0.98}}
\multiput(69.49,0.19)(0.97,0.06){1}{\line(1,0){0.97}}
\multiput(70.46,0.25)(0.96,0.07){1}{\line(1,0){0.96}}
\multiput(71.42,0.32)(0.95,0.08){1}{\line(1,0){0.95}}
\multiput(72.38,0.40)(0.95,0.08){1}{\line(1,0){0.95}}
\multiput(73.32,0.48)(0.94,0.09){1}{\line(1,0){0.94}}
\multiput(74.26,0.57)(0.92,0.10){1}{\line(1,0){0.92}}
\multiput(75.18,0.67)(0.91,0.11){1}{\line(1,0){0.91}}
\multiput(76.10,0.78)(0.90,0.12){1}{\line(1,0){0.90}}
\multiput(77.00,0.90)(0.89,0.12){1}{\line(1,0){0.89}}
\multiput(77.88,1.02)(0.87,0.13){1}{\line(1,0){0.87}}
\multiput(78.76,1.15)(0.86,0.14){1}{\line(1,0){0.86}}
\multiput(79.61,1.29)(0.84,0.14){1}{\line(1,0){0.84}}
\multiput(80.45,1.43)(0.82,0.15){1}{\line(1,0){0.82}}
\multiput(81.28,1.58)(0.81,0.16){1}{\line(1,0){0.81}}
\multiput(82.08,1.74)(0.79,0.17){1}{\line(1,0){0.79}}
\multiput(82.87,1.91)(0.77,0.17){1}{\line(1,0){0.77}}
\multiput(83.64,2.08)(0.75,0.18){1}{\line(1,0){0.75}}
\multiput(84.39,2.26)(0.36,0.09){2}{\line(1,0){0.36}}
\multiput(85.12,2.44)(0.35,0.10){2}{\line(1,0){0.35}}
\multiput(85.82,2.63)(0.34,0.10){2}{\line(1,0){0.34}}
\multiput(86.51,2.83)(0.33,0.10){2}{\line(1,0){0.33}}
\multiput(87.17,3.03)(0.32,0.10){2}{\line(1,0){0.32}}
\multiput(87.81,3.24)(0.31,0.11){2}{\line(1,0){0.31}}
\multiput(88.42,3.45)(0.30,0.11){2}{\line(1,0){0.30}}
\multiput(89.02,3.67)(0.28,0.11){2}{\line(1,0){0.28}}
\multiput(89.58,3.89)(0.27,0.11){2}{\line(1,0){0.27}}
\multiput(90.12,4.12)(0.26,0.12){2}{\line(1,0){0.26}}
\multiput(90.63,4.35)(0.24,0.12){2}{\line(1,0){0.24}}
\multiput(91.12,4.59)(0.23,0.12){2}{\line(1,0){0.23}}
\multiput(91.58,4.83)(0.22,0.12){2}{\line(1,0){0.22}}
\multiput(92.01,5.07)(0.20,0.12){2}{\line(1,0){0.20}}
\multiput(92.42,5.32)(0.19,0.13){2}{\line(1,0){0.19}}
\multiput(92.79,5.57)(0.17,0.13){2}{\line(1,0){0.17}}
\multiput(93.14,5.82)(0.16,0.13){2}{\line(1,0){0.16}}
\multiput(93.46,6.08)(0.14,0.13){2}{\line(1,0){0.14}}
\multiput(93.75,6.34)(0.13,0.13){2}{\line(0,1){0.13}}
\multiput(94.01,6.60)(0.11,0.13){2}{\line(0,1){0.13}}
\multiput(94.24,6.86)(0.10,0.13){2}{\line(0,1){0.13}}
\multiput(94.44,7.13)(0.17,0.27){1}{\line(0,1){0.27}}
\multiput(94.61,7.39)(0.14,0.27){1}{\line(0,1){0.27}}
\multiput(94.75,7.66)(0.11,0.27){1}{\line(0,1){0.27}}
\multiput(94.86,7.93)(0.08,0.27){1}{\line(0,1){0.27}}

\thicklines
\path(95.00,8.13)(95.63,6.25)

\thicklines
\path(95.00,8.13)(93.13,7.50)

\end{picture}
 \caption{$z_i$ goes around $z_j$ and $z_k$ by counter-clockwise}
 \label{fig:i(jk)}
\end{figure}

\begin{proof}
Proof of \eqref{CON}. \ We have
\begin{align*}
 W_{(ij)k}(z_1, z_2, z_3)
 &= G_0\left(X_{ij},-X_{jk};\frac{z_{ji}}{z_{ki}}\right)
    \times z_{ki}^{T}, \\
 W_{i(jk)}(z_1, z_2, z_3)
 &= G_1\left(X_{ij},-X_{jk};\frac{z_{ji}}{z_{ki}}\right)
    \times z_{ki}^T %\\
\end{align*}
Hence we have
\begin{align*}
 W_{(ij)k}(z_1, z_2, z_3)
 & =  G_0\left(X_{ij},-X_{jk};\frac{z_{ji}}{z_{ki}}\right) \times z_{kj}^T\\
 & = G_1\left(X_{ij},-X_{jk};\frac{z_{ji}}{z_{ki}}\right)
            \times \varphi_{\mathrm{KZ}}(X_{ij}, -X_{jk})
                             \times z_{kj}^T\\
 & =  W_{i(jk)} \times \varphi_{\mathrm{KZ}}(X_{ij},-X_{jk}).
\end{align*}
Thus we have \eqref{CON}.

\noindent
Proof of \eqref{CONNE}. \ We have
\begin{align*}
 & W_{(ji)k}(z_1, z_2, z_3)\\
 & = G_0\left( X_{ji},-X_{ik};\frac{z_{ij}}{z_{kj}} \right) \times z_{kj}^T\\
 & = G_0\left(X_{ij},-X_{ik};\frac{z_{ji}/z_{ki}}{z_{ji}/z_{ki}-1} \right)
    \times z_{kj}^T\\
 & =\left(\frac{z}{z-1}\right)_* \circ \pi\left(\frac{z}{z-1}\right)
       \left(G_0\left(X_{ij},-X_{ik};\frac{z_{ji}}{z_{ki}}\right)\right)
    \times z_{kj}^T\\
 & =\left(\frac{z}{z-1}\right)_*
 \left(
  G_0\left( X_{ij},-X_{ik};\frac{z_{ji}}{z_{ki}} \right) \times \exp(\mp X_{ij} \pi i)
 \right) \times z_{kj}^T\\
 & =
  G_0\left( X_{ij},X_{ij}+X_{ik};\frac{z_{ji}}{z_{ki}} \right)
   \times \exp(\mp X_{ij} \pi i)
    \times z_{kj}^T\\
 & =
  \left( \frac{z_{kj}}{z_{ki}} \right)^{-X_{ij}-X_{ik}}
   \overline{G}_0\left( X_{ij},X_{ij}+X_{ik};\frac{z_{ji}}{z_{ki}} \right)
    \left( \frac{z_{ji}}{z_{kj}}\right)^{X_{ij}} 
   \times \exp(\mp X_{ij}\pi i)
    \times z_{kj}^T.
\end{align*}

Because $\sum_W wW$ is exponential of Lie series,
$\overline{G}_0$ is also exponential of Lie series w.r.t.\ $X$ and $Y$.
So adding any central element of Lie algebra to $X$ or $Y$
does not change $\overline{G}_0$, in particular
\begin{align*}
  \overline{G}_0\left( X_{ij},X_{ij}+X_{ik};\frac{z_{ji}}{z_{ki}} \right)
 =
  \overline{G}_0\left( X_{ij},X_{ij}+X_{ik}-T;\frac{z_{ji}}{z_{ki}} \right)
 =
  \overline{G}_0\left( X_{ij},-X_{jk};\frac{z_{ji}}{z_{ki}} \right).
\end{align*}
Using this equation we obtain
\begin{align*}
& W_{(ji)k}(z_1, z_2, z_3)\\
 & =
  \left( \frac{z_{kj}}{z_{ki}} \right)^{X_{jk}}
   \overline{G}_0\left( X_{ij},-X_{jk};\frac{z_{ji}}{z_{ki}} \right)
    \left(\frac{z_{ji}}{z_{kj}}\right)^{X_{ij}}
     \times z_{ki}^T
      \times \exp(\mp X_{ij} \pi i)\\
 & =
  W_{(ij)k}(z_1,z_2,z_3)   \times \exp(\mp X_{ij} \pi i),
\end{align*}
where $\frac{z_{ji}}{z_{ki}} \in \mathbf{H}_{\pm}$.
Thus we have \eqref{CONNE}.

\noindent
Proof of \eqref{CONNEC}. \ We have
\begin{align*}
 &W_{i(kj)}(z_1, z_2, z_3)\\
 &= G_1\left( X_{ik},-X_{kj}; \frac{z_{ki}}{z_{ji}} \right) \times z_{ji}^T\\
 &= G_0\left( X_{ik},-X_{kj}; \frac{1}{z_{ji}/z_{ki}} \right)
     \times \varphi_{\mathrm{KZ}}(X_{ik},-X_{kj})^{-1}
      \times z_{ji}^T\\
 &= \left(\frac{1}{z}\right)_* \circ \pi\left(\frac{1}{z}\right)
       \left(G_0\left(X_{ik},-X_{kj};\frac{z_{ji}}{z_{ki}}\right)\right)
      \times \varphi_{\mathrm{KZ}}(X_{ik},-X_{kj})^{-1}
       \times z_{ji}^T\\
 \begin{split}
 &= 
    \left.
    \left(\frac{1}{z}\right)_*
     \left(
      G_0\left(X,Y;\frac{z_{ji}}{z_{ki}}\right)
       \varphi_{\mathrm{KZ}}(X,Y)^{-1}
        \exp(\mp Y\pi i)
         \varphi_{\mathrm{KZ}}(-Y,X-Y)^{-1}
     \right)
   \right|_{\begin{subarray}{l}
    X=X_{ik}\\ Y=-X_{kj}
	    \end{subarray}}\\
 &\qquad
     \times \varphi_{\mathrm{KZ}}(X_{ik},-X_{kj})^{-1}
      \times z_{ji}^T
 \end{split}\\
 \begin{split}
 &= 
  \left.
    G_0\left(-X+Y,Y;\frac{z_{ji}}{z_{ki}}\right)
       \varphi_{\mathrm{KZ}}(-X+Y,Y)^{-1}
        \exp(\mp Y\pi i)
         \varphi_{\mathrm{KZ}}(-Y,-X)^{-1}
   \right|_{\begin{subarray}{l}
    X=X_{ik}\\ Y=-X_{kj}
	    \end{subarray}}\\
 &\qquad
     \times \varphi_{\mathrm{KZ}}(X_{ik},-X_{kj})^{-1}
      \times z_{ji}^T
 \end{split}\\
 &= 
    G_0\left(-X_{ik}-X_{kj},-X_{kj};\frac{z_{ji}}{z_{ki}}\right)
       \varphi_{\mathrm{KZ}}(X_{ij},-X_{kj})^{-1}
       \times z_{ji}^T
        \times \exp(\pm X_{kj} \pi i)
      \\
\intertext{(in the same way as above)}
 &= 
    G_0\left(X_{ij},-X_{kj};\frac{z_{ji}}{z_{ki}}\right)
       \varphi_{\mathrm{KZ}}(X_{ij},-X_{kj})^{-1}
       \times z_{ki}^T
        \times \exp(\pm X_{kj} \pi i)
      \\
 &= 
    G_1\left(X_{ij},-X_{kj};\frac{z_{ji}}{z_{ki}}\right)
     \times z_{ki}^T \times \exp(\pm X_{kj} \pi i) \\
 &=
  W_{i(jk)}(z_1, z_2, z_3)
   \times \exp(\pm X_{jk} \pi i)
\end{align*}
where $\frac{z_{ji}}{z_{ki}}\in\mathbf{H}_{\pm}$.
Thus we have \eqref{CONNEC}

Proof of \eqref{CONNECT}. \ We have
\begin{align*}
 &W_{(ij)k}(z_1,z_2,z_3)\\
 &=
 G_0 \left(X_{ij},-X_{jk};\frac{z_{ji}}{z_{ki}}\right)
  \times z_{ki}^T\\
 &=
 G_0
  \left(X_{ij},-X_{jk};
   \frac{\frac{z_{ik}}{z_{jk}}-1}{\frac{z_{ik}}{z_{jk}}}
  \right)
  \times z_{ki}^T\\
 &=
 \left(\frac{z-1}{z}\right)_*
 \circ
 \pi\left(\frac{1}{1-z}\right)
 \left(
  G_0
   \left(X_{ij},-X_{jk};
    \frac{z_{ik}}{z_{jk}}
   \right)
  \right)
  \times z_{ki}^T\\
 &=
 \left(\frac{z-1}{z}\right)_*
 \left(
  G_1
   \left(X_{ij},-X_{jk};
    \frac{z_{ik}}{z_{jk}}
   \right)\times \exp(\pm X_{jk}\pi i)
  \right)
  \times z_{ki}^T\\
 &=
  G_1
   \left(-X_{ij}-X_{jk},-X_{ij};
    \frac{z_{ik}}{z_{jk}}
   \right)\times \exp(\pm X_{ij}\pi i)
  \times z_{ki}^T\\
 \begin{split}
  &=
  \left(\frac{z_{ji}}{z_{jk}}\right)^{X_{ij}}
  \overline{G}_1
   \left(X_{ki},-X_{ij};
    \frac{z_{ik}}{z_{jk}}
   \right)
  \left(\frac{z_{ji}}{z_{jk}}\right)^{-X_{ij}-X_{jk}}\\
  &\hspace{4cm}
  \times \exp(\pm X_{ij}\pi i)
  \times \left(\frac{z_{ki}}{z_{jk}}\right)^T
  \times z_{jk}^T.
 \end{split}
\end{align*}
Here we introduce $\overline{G}_1$ by
\begin{align*}
  G_1(X,Y;z) = z^X \times \overline{G}_1(X,Y;z) \times (1-z)^{-Y}.
\end{align*}
$\overline{G}_1$ is also exponential of Lie series,
in the same way above, we get
\begin{align*}
&W_{(ij)k}(z_1,z_2,z_3)\\
 &=
  \left(\frac{z_{ji}}{z_{jk}}\right)^{X_{ij}}
  \overline{G}_1
   \left(X_{ki},-X_{ij};
    \frac{z_{ik}}{z_{jk}}
   \right)
  \left(\frac{z_{ji}}{z_{jk}}\right)^{X_{ik}}
  \times z_{jk}^T
  \times \exp(\mp (X_{ik}+X_{jk})\pi i)\\
 &=
 W_{k(ij)}(z_1,z_2,z_3)  \times \exp(\mp (X_{ik}+X_{jk})\pi i).
\end{align*} 
where $\frac{z_{ji}}{z_{ki}}\in\mathbf{H}_{\pm}$.
Thus we have \eqref{CONNECT}.

The last equation in the proposition can be shown as follows:
\begin{align*}
 &W_{i(jk)}(z_1, z_2, z_3)\\
 &=
  G_1\left(X_{ij},-X_{jk};\frac{z_{ji}}{z_{ki}}\right)
   \times z_{ki}^T\\
 &=
  G_0\left(X_{ij},-X_{jk};\frac{z_{ji}}{z_{ki}}\right)
   \times \varphi_{\mathrm KZ}(X_{ij},-X_{jk})^{-1}
   \times z_{ki}^T\\
 &=
  \left(\frac{1}{1-z}\right)_*
   \circ \pi\left(\frac{z-1}{z}\right)
    \left(
     G_0\left(X_{ij},-X_{jk};\frac{z_{kj}}{z_{ij}}\right)
    \right)
     \times \varphi_{\mathrm KZ}(X_{ij},-X_{jk})^{-1}
      \times z_{ki}^T\\
 &=
  \left(\frac{1}{1-z}\right)_*
    \left(
     G_\infty\left(X_{ij},-X_{jk};\frac{z_{kj}}{z_{ij}}\right)
      \times \exp(\pm (-X_{ij}-X_{jk})\pi i)
    \right)\\
 &\hspace{7cm}
  \times \varphi_{\mathrm KZ}(X_{ij},-X_{jk})^{-1}
      \times z_{ki}^T\\
 &=
   G_\infty\left(X_{jk},X_{ij}+X_{jk};\frac{z_{kj}}{z_{ij}}\right)
      \times \exp(\pm X_{ij}\pi i)
  \times \varphi_{\mathrm KZ}(X_{ij},-X_{jk})^{-1}
      \times z_{ki}^T\\
 &=
  G_0\left(X_{jk},X_{ij}+X_{jk};\frac{z_{kj}}{z_{ij}}\right)
   \exp(\mp X_{jk}\pi i)
    \varphi_{\mathrm{KZ}}(X_{jk},-X_{ij})^{-1}
     \exp(\mp X_{ij}\pi i)\\
 &\hspace{2cm}
     \times \exp(\pm X_{ij}\pi i)
    \times \varphi_{\mathrm KZ}(X_{ij},-X_{jk})^{-1}
      \times z_{ki}^T\\
 &=
  G_0\left(X_{jk},X_{ij}+X_{jk};\frac{z_{kj}}{z_{ij}}\right)
   \times \exp(\mp X_{jk}\pi i)
      \times z_{ki}^T\\
 &=
  G_0\left(X_{jk},-X_{ki};\frac{z_{kj}}{z_{ij}}\right)
   \times z_{ij}^T
    \times \exp(\mp (-X_{ij}-X_{ik})\pi i)\\
 &=
  W_{(jk)i}(z_1, z_2, z_3)
    \times \exp(\pm (X_{ij}+X_{ik}))\pi i),
\end{align*}
where $\frac{z_{ji}}{z_{ki}}\in \mathbf{H}_{\pm}$.

\end{proof}

From this proposition,
by comparing two ways to analytic continuation as in \cite{Kas},
we obtain so called hexagon relation for the Drinfel$'$d associator

By substituting $\varphi$ satisfy
 \begin{align*}
  &\exp( (X_{ik} + X_{jk})\pi i)\\
  &\quad
  = \varphi_{\mathrm{KZ}}(X_{ik},X_{ij}) \exp(X_{ik}\pi i) \varphi_{\mathrm{KZ}}(X_{ik},X_{jk})^{-1}
    \exp(X_{jk}\pi i) \varphi_{\mathrm{KZ}}(X_{ij},X_{jk}),\\
  &\exp((X_{ij}+X_{ik})\pi i)\\
  &\quad
  = \varphi_{\mathrm{KZ}}(X_{jk},X_{ik})^{-1} \exp(X_{ik}\pi i) \varphi_{\mathrm{KZ}}(X_{ij},X_{ik})
  \exp(X_{ij}\pi i) \varphi_{\mathrm{KZ}}(X_{ij},X_{jk})^{-1}.
 \end{align*}
Thus using Landen and Euler connection formulas,
we obtain the hexagon relations.

\section*{Appendix 2: 
Connection formulas for the confluent hypergeometric functions and the functional
relation for the Hurwitz zeta function}

\begin{center}
Kimio Ueno and Michitomo Nishizawa
\end{center}

\vspace{1cm}

The generalized zeta function introduced by Hurwitz (Hurwitz zeta function, for short)
is, by definition, 
\begin{eqnarray}
   \zeta(s,z) \,=\, \sum_{k=0}^{\infty}\frac{1}{(k+z)^s},
\end{eqnarray}
where we suppose that $0 \leq z < 1$. The series is absolutely convergent for $\Re s > 1$,
and is analytically continued to the whole $s$-plane as a meromorphic function. 
Evidently, $\zeta(s,1) \,=\, \zeta(s)$, which is the Riemann Zeta function. 
Furthermore, they satisfy a functional relation
\begin{eqnarray}\label{HURWITZ}
    \zeta(s,z) \,=\, \Gamma(1-s) \left\{(2\pi i)^{s-1}\mathcal{L}(1-s,z)
                      +(-2\pi i)^{s-1}\mathcal{L}(1-s,1-z) \right\},
\end{eqnarray}
which was established by Hurwitz himself \cite{WW}. Here 
\begin{eqnarray*}
      \mathcal{L}(s,z) \,=\, \sum_{n=1}^{\infty}\frac{e^{2\pi inz}}{n^s}
\end{eqnarray*}
is a generalized polylogarithm and $\Gamma(s)$ is Euler's gamma function.
Let call (\ref{HURWITZ}) the Hurwitz relation.
When $z=1$, it reduces to the functional equation for the Riemann zeta function:
\begin{eqnarray*}
  \zeta(s) \,=\, 2^s\pi^{s-1}\Gamma(1-s)\sin\left(\frac{\pi s}{2}\right)\zeta(1-s).
\end{eqnarray*}
To investigate analytic continuation of zeta functions, the sum formula of Euler-Maclaurin
\begin{eqnarray}\label{EM}
     \sum_{r=0}^{\infty}f(r) \,=\, \int_0^{\infty}\!\!f(t)\,dt
       +\sum_{k=1}^n \, \frac{B_k}{k!}\{ f^{(k-1)}(\infty)-f^{(k-1)}(0) \} \nonumber \\
       +\frac{(-1)^n}{n!}\int_0^{\infty}\!\!\overline{B}_n(t)f^{(n)}(t)\,dt
\end{eqnarray}
is a useful tool. Here $\overline{B}_n(t)=B_n(t-[t])$ ($[t]$ denotes the integral 
part of $t$),
$B_n(t)$ is the n-th Bernoulli polynomial defined by
\begin{eqnarray*}
         \sum_{n=0}^{\infty}\frac{B_n(t)u^n}{n!} \,=\,
	  %\sum_{n=0}^{\infty}
	  \frac{{}ue^{tu}{}}{e^u-1},
\end{eqnarray*}
and $B_n=B_n(0)$ is the n-th Bernoulli number. Putting $f(t)=(t+z)^{-s}$, 
and $n=2$ in (\ref{EM}), we have
\begin{eqnarray}\label{EM1}
  \zeta(s,z) \,=\, \frac{z^{1-s}}{s-1}+\frac{z^{-s}}{2}+\frac{sz^{-s-1}}{12}
            - \int_0^{\infty}\frac{\overline{B}_2(t)}{2!}
                \left(\frac{d}{dt}\right)^2\left\{\frac{1}{(z+t)^s}\right\}\,dt
\end{eqnarray}
Substituting the Fourier expansion
\begin{eqnarray}
       \overline{B}_2(t)=\frac{1}{\pi^2}\sum_{n=0}^{\infty}\frac{\cos(2\pi nt)}{n^2}
\end{eqnarray}
into (\ref{EM1}), making once partial integration, we obtain
\begin{eqnarray}\label{AN}
  \zeta(s,z) \,=\, \frac{z^{1-s}}{s-1}+\frac{z^{-s}}{2}+\frac{sz^{-s}}{2\pi i}
        \sum_{l\neq0}\frac{1}{\,l\,}
                 \int_0^{\infty}\frac{e^{-2\pi il zu}}{(1+u)^{s+1}}\,du.
\end{eqnarray}

Now we observe the confluent hypergeometric equation
\begin{eqnarray}\label{CONFL}
   x\frac{d^2y}{dx^2}+(\gamma-x)\frac{dy}{dx}-\alpha y = 0.
\end{eqnarray}
This equation has a regular singular point at $x=0$, and an irregular singular point
at $x=\infty$. The confluent hypergeometric series
\begin{eqnarray}\label{CONFLFUNC}
    F(\alpha,\gamma;x) \,=\, \sum_{n=0}^{\infty} 
             \frac{(\alpha)_nx^n}{(\gamma)_nn!},
\end{eqnarray}
and \ $x^{1-\gamma}F(\alpha-\gamma+1,2-\gamma;x) \, (\,=e^xx^{1-\gamma}F(1-\alpha,2-\gamma;-x))$
\ form a system of fundamental solutions to (\ref{CONFL}) around $x=0$.
Let us introduce a function defined by
\begin{eqnarray}\label{U}
  U(\alpha,\gamma;x) \,=\, \frac{1}{\Gamma(\alpha)}
                     \int_0^{\infty}e^{-xu}(1+u)^{\gamma-\alpha-1}u^{\alpha-1}\,du.
\end{eqnarray}
This is a solution to (\ref{CONFL}) around $x=\infty$ \cite{S}, 
and is connected to the former solutions by
\begin{eqnarray}\label{CONN}
     U(\alpha,\gamma;x) \,=\, \frac{\Gamma(1-\gamma)}{\Gamma(1+\alpha)}F(\alpha,\gamma;x)
         + \frac{\Gamma(\gamma-1)}{\Gamma(\alpha)}e^xx^{1-\gamma}F(1-\alpha,2-\gamma;-x)).
\end{eqnarray}
Furthermore, for $\Re \alpha>0$, we see that
\begin{eqnarray}\label{ASYM}
        U(\alpha,\gamma;x) \,\sim\, x^{-\alpha}, \qquad (|x| \to \infty).
\end{eqnarray}
Substituting (\ref{U}) to (\ref{AN}) with $\alpha=1, \gamma=1-s, \text{and} \ x=-2\pi il z$,
we obtain
\begin{eqnarray}\label{EM2}
  \zeta(s,z) \,=\, \frac{z^{1-s}}{s-1}+\frac{z^{-s}}{2}+\frac{sz^{-s}}{2\pi i}
        \sum_{l\neq0}\frac{1}{\,l\,}U(1,1-s;-2\pi i l z).
\end{eqnarray}
Here we should note that the infinite sum in (\ref{EM2}) is absolutely convergent due to 
the asymptotic behavior (\ref{ASYM}). According to (\ref{CONN}), 
we have a connection formula
\begin{eqnarray}
U(1,1-s;-2\pi il z) \,=\, \frac{1}{s}\left\{ F(1,1-s;-2\pi il z) - 
                              \Gamma(1-s)(-2\pi il z)^se^{-2\pi il z} \right\}.
\end{eqnarray}
Hence
\begin{eqnarray}
  \frac{sz^{-s}}{2\pi i}\sum_{l\neq0}\frac{1}{\,l\,}U(1,1-s;-2\pi i l z) = 
  \frac{z^{-s}}{2\pi i}\sum_{l\neq0}\frac{1}{\,l\,}F(1,1-s;-2\pi i l z)+ \nonumber \\
  + \Gamma(1-s) \left\{(2\pi i)^{s-1}\mathcal{L}(1-s,z)
                      +(-2\pi i)^{s-1}\mathcal{L}(1-s,1-z) \right\}.
\end{eqnarray}
By virtue of the integral representation for $F(\alpha,\gamma;x)$
\begin{eqnarray*}
  F(\alpha,\gamma;x) \,=\, \frac{\Gamma(\gamma)}{\Gamma(\gamma-\alpha)\Gamma(\alpha)}
        \int_0^1 e^{zt}t^{\alpha-1}(1-t)^{\gamma-\alpha-1}\,dt,
\end{eqnarray*}
we have
\begin{eqnarray*}
\frac{1}{2\pi i}\sum_{l\neq0}\frac{1}{\,l\,}F(1,1-s;-2\pi i l z) &=&
\frac{s}{\pi} \sum_{n=1}^{\infty} \int_0^1 \frac{\sin(2\pi nzt)}{n}(1-t)^{-s} \,dt \nonumber \\
&=& \frac{s}{\pi^2z} \left\{\frac{\pi^2}{6}+(s+1)
                \int_0^1 \overline{B}_2(t)(1-t)^{-s-2}\,dt \right\}.
\end{eqnarray*}
Since $\overline{B}_2(t)=\pi^2(t^2-t+1/6) \quad (0 \leq t \leq 1)$, we have
\begin{eqnarray*}
\frac{z^{1-s}}{s-1}+\frac{z^{-s}}{2}+\frac{z^{-s}}{2\pi i}
             \sum_{l\neq0}\frac{1}{\,l\,}F(1,1-s;-2\pi i l z) \,=\, 0
\end{eqnarray*}
for $0 \leq z <1$. Thus we obtain the Hurwitz relation (\ref{HURWITZ}).

In the case of the q-analogue of the Hurwitz zeta function, the connection formulas for the
hypergeometric function play the same role as in the present case \cite{UN1}.

\vspace{1cm}

\begin{tabular}{l}
 Jun-ichi Okuda and Kimio Ueno\\
 Department of Mathematical Sciences\\
 School of Science and Engineering\\
 Waseda University\\
 Okubo 3-4-1, Shinjuku-ku, Tokyo 169-8555, Japan\\
 E-mail: {\ttfamily okuda@gm.math.waseda.ac.jp}, {\ttfamily uenoki@mse.waseda.ac.jp}
\end{tabular} 

\vspace{1cm}

\begin{tabular}{l}
 Michitomo Nishizawa\\
 Graduate School of Mathematical Sciences\\
 University of Tokyo\\
 Komaba 3-8-1, Meguro-ku, Tokyo 153-8914, Japan\\
 E-mail: {\ttfamily mnishi@ms.u-tokyo.ac.jp}
\end{tabular}

\end{document}